%% file: paperdraft_07.tex
\documentclass{siamart250211}

\usepackage[normalem]{ulem}
\usepackage{addedsty}

\DefineNamedColor{named}{JungleGreen} {cmyk}{0.99,0,0.52,0}

\input{ex_shared}

\usepackage{ifpdf}
\ifpdf
\hypersetup{
  pdftitle={Christoffel Adaptive Sampling for Sparse Random Features Expansion},
  pdfauthor={}
}
\fi

\begin{document}

\maketitle

\begin{abstract}
Random Feature Models (RFMs) have become a powerful tool for approximating multivariate functions and solving partial differential equations efficiently. Sparse Random Feature Expansions (SRFE) improve traditional RFMs by incorporating sparsity, making it particularly effective in data-scarce settings. In this work, we integrate active learning with sparse random feature approximations to improve sampling efficiency. Specifically, we incorporate the Christoffel function to guide an adaptive sampling process, dynamically selecting informative sample points based on their contribution to the function space. This approach optimizes the distribution of sample points by leveraging the Christoffel function associated with an iteratively-chosen basis obtained by the sparse recovery solver. We conduct numerical experiments comparing adaptive and non-adaptive sampling strategies with the SRFE framework and examine their accuracy for various function approximation tasks. Overall, our results demonstrate the advantages of adaptive sampling in maintaining high accuracy while reducing sample complexity for SRFE, highlighting its potential for scientific computing tasks where data is expensive to acquire.
\end{abstract}



\section{Introduction}
Many scientific computing processes require approximating complex functions in more than one dimension, where obtaining data is difficult and uniform sampling quickly becomes inefficient. We provide an adaptive approach to improve the sampling process. In contrast to the standard machine learning setting, which draws a single fixed dataset, this approach allows the training data to be selected strategically and then reused to refine the sampling strategy. Such a setting is reasonable in many applications where data acquisition is expensive or limited, such as numerical simulations, physical experiments, or scientific computing tasks that require function evaluations at high computational cost.

\subsection{Overview}
The goal of this work is to improve multivariate function approximation through an \textit{adaptive sampling} process guided by a designed selection criterion. Here adaptive sampling refers to the process in which samples are selected dynamically across iterations: existing samples are reused, and new samples are chosen based on the function approximation obtained in the previous iteration. 

Throughout this work, we consider a domain $D \subset \mathbb{R}^d, d \geq 1$, with underlying probability measure $\rho$, and an unknown target function $f: D \rightarrow \mathbb{R}$. Given numbers $1 \leq m_1 < m_2 < \cdots$, we approximate the function using a nested collection of sample points (previous samples are retained across iterations):
\begin{equation}
\{x_k\}_{k=1}^{m_1} \subseteq \{x_k\}_{k=1}^{m_2} \subseteq \cdots
\label{eq:adapsamples}
\end{equation}
Using the training data $\{(x_k, f(x_k))\}_{k=1}^{m_i}$, $i = 1,2,\ldots$, we compute a sequence of approximations $\hat{f}_1,\hat{f}_2,\ldots$. The procedure we implement involves both \textit{adaptive approximation} and \textit{adaptive sampling}. At step $i$, given the approximation $\hat{f}_i$ is used to compute new sample points $\{ x_k \}^{m_{i+1}}_{k=1}$ (adaptive sampling), which are then used to compute the new approximation $\hat{f}_{i+1}$ (adaptive approximation). In doing so, we aim to achieve rapid decrease of the approximation error $\|f - \hat{f_i}\|_{L^2_\rho(D)}$ as the number of iterations $i \rightarrow \infty$. In particular, we aim for faster decrease over standard \textit{nonadaptive} sample procedures -- specifically, \textit{Monte Carlo (MC)} sampling, where sample points $x_k$ are drawn i.i.d.\ from the underlying probability measure $\rho$.

\subsection{Main contributions}
In scientific machine learning, data scarcity arises in many applications, such as Physics-Informed Neural Networks (PINNs) for PDEs \cite{han2018solving, lu2021learning, raissi2019physics-informed, sirignano2018dgm}, operator learning \cite{bhattacharya2021model, li2021fourier, 1964handbook, gajjar2022provable, geist2021numerical, heiss2021neural, kovachki2023neural, kutyniok2022theoretical, wang2021learning}, image recovery \cite{adcock2021compressive, liang2020deep, mccann2019biomedical, ongie2020deep, ravishankar2020image, sandino2020compressed} and physical system modeling \cite{rudy2017data-driven, brunton2024data-driven, brunton2023machine}. It is therefore imperative to design highly sample-efficient algorithms.

To cope with the issue of data scarcity, we implement a representation framework that is effective in approximating functions in the low-data regime. The Random Feature Model (RFM) is a well-known and increasingly widely-used technique approximating multivariate functions in scientific machine learning. However, for high accuracy, RFM usually requires a large dataset -- in particular, more data than trainable parameters. Sparse Random Feature Expansion (SRFE) was introduced in \cite{hashemi2023generalization} to address the issue of limited data by incorporating sparsity to improve data efficiency. However, prior work on SRFE considers only nonadaptive, MC sampling. 

Our main contribution in this work is to enhance the accuracy of SRFE in the limited data regime through adaptive sampling. Our procedure is based on Christoffel Adaptive Sampling (CAS), which is an extension of so-called Christoffel Sampling (CS) to nonlinear approximation spaces. CS is an active learning technique for (weighted) least-squares approximation in linear approximation spaces. While CS provides near-optimal recovery guarantees for linear spaces, many machine learning tasks involve nonlinear spaces. CAS facilitates the use of CS for nonlinear spaces, by dynamically implementing CS as part of an adaptive sampling process.

We combine CAS and SRFE in this work, introducing Christoffel Adaptive Sampling for Sparse Random Feature Expansions (CAS-SRFE). A major challenge we address is efficiently computing the samples generated by CAS. Differing from previous work, which uses inefficient constructions based on finite grids, we implement a novel method which combines a suitable Gram matrix and an efficient sampler based on the Metropolis--Hastings (MH) algorithm. Finally, we present a series of numerical experiments to showcase the efficiency of CAS-SRFE, demonstrating throughout a consistent benefit of adaptive sampling over that of nonadaptive, MC sampling. Overall, we demonstrate how the CAS-SRFE framework unites sparsity, adaptivity and active learning to effectively address data-scarce problems in scientific computing.

\subsection{Relevant literature}

The foundation for this paper is the well-known RFM. RFM uses simple functions with randomized parameters to approximate functions, and can be viewed as an ensemble average of randomly parameterized functions \cite{nelsen2024operator}. 
Rahimi and Recht \cite{rahimi2008random} originally introduced RFM as an efficient approximation to kernel methods, offering better scalability with respect to the size of the training dataset \cite{rahimi2008random}. It was later shown to be successful for many machine learning tasks, with the benefits of numerical efficiency and theoretical accuracy bounds \cite{hashemi2023generalization}. For instance, this approximation is also used downstream for kernel regression tasks \cite{hashemi2023generalization}. The most common form of RFM is the Random Fourier Feature (RFF) method, which uses trigonometric functions as features. In our numerical experiments, we only consider RFF, although CAS-SRFE can be readily applied to other features.

While effective, RFM relies on large datasets, as it requires more samples than trainable parameters. SRFE was introduced by Hashemi et al.\ \cite{hashemi2023generalization} to address this limitation. SRFE incorporates sparse recovery into RFM. It constructs a sparse representation using a large set of candidate features from a given set. By doing so, one can allow the number of features to significantly exceed the dataset size. This provides better representation capabilities than RFM and reasonable generalization bounds, even in the highly data-scarce setting.

The original SRFE framework solved a convex $\ell^1$-minimization problem, which is known not to yield sparse approximations in general \cite[Ex.\ 3.2]{foucart2013mathematical}. Since our CAS-SRFE framework relies on obtaining a sparse approximation at each step, we utilize Hard Thresholding Pursuit (HTP) and Orthogonal Matching Pursuit (OMP) instead  \cite[Chpt.\ 3]{foucart2013mathematical}. Besides yielding sparse solutions, these are simpler, as they do not require a convex optimization solver. They also may be more computationally in practice, although that is not the primary concern of this work.

SRFE has been extended in various ways. In \cite{saha2023harfe}, the authors proposed Hard-Ridge Random Feature Expansion (HARFE), a sparse ridge regression model solved via modified version of the HTP algorithm, and compared its performance with SRFE. In \cite{xie2022shrimp}, the authors introduced Sparser Random Feature Models via Iterative Magnitude Pruning (ShRIMP) to more efficiently fit high-dimensional data with low-dimensional structure. The algorithm in HARFE is similar to our CAS-SRFE with HTP as the sparse recovery algorithm, albeit without adaptive sampling. Both HARFE and ShRIMP could be deployed within our CAS-SRFE framework. However, since our focus in this work is to explore how to incorporate adaptive sampling into SRFE, we focus on simpler underlying sparse recovery algorithms. We leave such extensions to future work (see \S\ref{sec:conclusions}).

The Christoffel function of a finite-dimensional, linear subspace is a classical object in approximation theory \cite{hampton2015coherence,nevai1986geza, xu1995christoffel}, that has recently become a powerful tool for sampling \cite{dolbeault2022optimal, shin2016near, narayan2017christoffel, adcock2025optimala, adcock2022cas4dl, adcock2023cs4ml}. Given an $n$-dimensional approximation space $\cP$, CS involves i.i.d.\ random sampling according to a probability measure with density proportional to the Christoffel function. When combined with a (weighted) least-squares fit, it yields quasi-optimal approximations $\cP$ with near-optimal sample complexity: namely, the number of samples required is $\mathcal{O}(n \log(n))$. See \cite{adcock2025optimala} for a review.

CAS \cite{adcock2022cas4dl} is a means to extend CS to nonlinear approximation. CAS dynamically updating samples via the Christoffel function of a judiciously- and adaptively-chosen sequence of linear approximation spaces.  It re-approximates the target function at each iteration, constructs the corresponding subspaces, and then updates the Christoffel function accordingly. 
CAS is both an adaptive approximation and an adaptive sampling method. This property is particularly advantageous when data is limited, as the adaptive procedure reuses samples to improve the approximation. CAS was originally developed for Deep Neural Network (DNN) approximation via Deep Learning (DL), and termed CAS4DL. In this work, we solidify CAS as a general framework for incorporating adaptive sampling into nonlinear approximation by demonstrating that it can also be combined with other popular nonlinear approximation schemes, namely, SRFE. It is worth stressing that while both use the same underlying procedure, i.e., CAS, CAS4DL and CAS-SRFE are very different in how they generate the requisite subspaces. CAS4DL uses a basis defined by the final layer of the DNN, while SRFE uses the basis corresponding to the selected sparse random features.

A major computational challenge in CS, and, consequently, CAS, is efficiently drawing samples from the resulting probability measure (the \textit{CS measure}, as we henceforth term it) \cite{adcock2022cas4dl,adcock2025optimala}.  The Christoffel function can be expressed explicitly in terms of any orthonormal basis of $\cP$. However, in practice, many approximation schemes, especially those based on DNN and RFM, do not yield explicit orthonormal bases, but rather, simply spanning sets which may also be redundant.

A common strategy in CS, also used in CAS, involves constructing finite grids as approximations of the domain.
As introduced in \cite{adcock2020near-optimal, migliorati2021multivariate}, constructing finite grids involves replacing the continuous domain $D$ with a finite but sufficiently fine grid, denoted by $Z = \{z_k\}_{k=1}^{K} \subset D$. An orthonormal basis is then constructed via linear algebra (e.g., QR factorization), operating on a dense $K \times n$ matrix, where $n = \mathrm{dim}(\cP)$. The resulting CS measure is then a discrete measure supported on $Z$. Crucially, this approach requires selecting a sufficiently large $K$ to ensure an accurate approximation off the grid, as further discussed in \cite{adcock2020near-optimal, migliorati2021multivariate, trunschke2024optimal,herremans2025refinement}. However, this can lead to prohibitive computational costs, as $K$ often needs to be exceedingly large.

In this work, we exploit the structure of RFM to work with the true CS measure, rather than a discretized variant, thereby avoiding the bottleneck arising from discretization. We first efficiently compute the Christoffel function by computing an eigenvalue decomposition of an $n \times n$ Gram matrix. We then adopt the MH algorithm for the sampling process. MH is widely used for its simplicity and versatility, providing a straightforward approach to such sampling tasks \cite{robert2016metropolis, robert1999monte}. In addition, compared to rejection sampling methods, the MH algorithm suffers less from the curse of dimensionality, as in the probability of rejection does not increase exponentially with the dimensions. It is therefore better suited to multivariate settings \cite{robert2016metropolis}.
As we demonstrate, this procedure allows for efficient implementation of CAS-SRFE.

\subsection{Outline}
In \S\ref{sec:background}, we provide background on RFM, SRFE and CS. In \S\ref{sec:main}, we describe the main algorithm, CAS-SRFE, including the implementation of each component. In \S\ref{sec:numerics}, we present a series of numerical experiments demonstrating the effectiveness of our framework, comparing the adaptive sampling process with the non-adaptive one. Finally,  \S\ref{sec:conclusions}, we conclude and discuss potential improvements. 
Code that implements CAS-SRFE is available here: \url{https://github.com/wangsherril/CAS4SRFE}.

\section{Background}
\label{sec:background}
We now describe RFM, SRFE and CS in more detail.

\subsection{RFM and SRFE}
\label{sec:rfm}

Following \cite{liao2026solving}, we consider a function $f: \mathbb{R}^d \rightarrow \mathbb{R}$ and, for the moment, fixed training samples  $\{(x_k, f(x_k))\}_{k=1}^{m}$. RFM approximates the function through a linear combination of a finite number of random features. We denote the random feature functions by $\phi(x;w)$, where $x \in \mathbb{R}^d$ and $w \in \mathbb{R}^d$. The total number of random features is denoted by $N$. To generate random features, we draw $w_1, \dots, w_N$ i.i.d. from a prescribed probability distribution $\gamma$.
We then construct an approximation $\hat{f}$ to the target function $f$ in the following form
\begin{equation}
    \hat{f}(\cdot) = \sum_{j=1}^N \hat{c}_j \phi(\cdot ; w_j).
    \label{RFM}
\end{equation}
The features $\{ \phi(\cdot; w_j) \}^{N}_{j=1}$ aim to form sufficiently rich families of functions to allow accurate recovery of the target function. Common examples for the feature function $\phi(x;w)$ include Random Fourier Features (RFF) $\phi({x}; {w}) = \exp(\mathrm{i}\langle {x}, {w} \rangle)$, random trigonometric features $\phi({x}; {w}) = \cos(\langle {x}, {w} \rangle)$ and random ReLU features $\phi({x}; {w}) = \max(\langle {x}, {w} \rangle, 0)$. As mentioned above, in this paper we focus on RFF.

Let $A$ be the random feature matrix and $b$ be the vector of function samples, i.e., 
\[
A= \frac{1}{\sqrt{m}} [\phi(x_s; w_r)]_{s,r=1}^{m, N} \in \mathbb{C}^{m \times N}, \quad 
{b}= \frac{1}{\sqrt{m}}[f(x_s)]^{m}_{s=1} \in \mathbb{C}^m.
\]
Training the RFM \eqref{RFM} is equivalent to finding the coefficient vector $\hat{c} \in \mathbb{R}^N$ such that $A \hat{c}\approx b$. In the under-parameterized regime ($m \geq N$), the coefficients are obtained by solving a (regularized) least squares problem. 
However, the focus of this paper is the over-parametrized regime ($N \geq m$), as data is often limited in scientific computing tasks. This is the setting of SRFE.
SRFE introduces the idea of sparsity from compressed sensing \cite{hashemi2023generalization} into RFM. 
In SRFE, one strives to construct a sparse coefficient vector satisfying $A \hat{c} \approx b$, for some sparsity level $s$ satisfying $1 \leq s \leq m \leq N$. This leads to a sparse recovery problem, which can be solved, for instance, by solving an $\ell^1$-minimization program. For example, \cite{hashemi2023generalization} considered the quadratically-constrained basis pursuit program
\[
\min_{c \in \mathbb{C}^N} \|c \|_{1}\text{ subject to } \| A c - b \|_{2} \leq \eta.
\]
Alternatively, one can use iterative and greedy methods, such OMP and HTP. Further details of the methods will be given in \S\ref{sec:step1}.

\subsection{Christoffel Sampling (CS)}
\label{sec:cs}
CS is a technique for constructing near-optimal sample points for (weighted) least-squares approximation in an arbitrary $n$-dimensional subspace.
Let $\mathcal{P} \subseteq L_{\rho}^2(D)$ with $\mathrm{dim}(\cP) = n$, where $\rho$ is the underlying probability measure on $D$. The Christoffel function $\mathcal{K}(x)$ of $\mathcal{P}$ is defined by
\begin{equation}
    \mathcal{K}(x) = \mathcal{K}(\mathcal{P})(x) := \sup \left\{ |p(x)|^2: p \in \mathcal{P}, \|p\|_{L_{\rho}^2(D)}^2 = 1\right\}, \quad \forall x \in D.
    \label{christoffeldef}
\end{equation}
Let \( \left\{\psi_j\right\}_{j=1}^{n} \) be any orthonormal basis of \( \mathcal{P} \).  Then $\cK(x)$ can be expressed as
\begin{equation}
    \mathcal{K}(x)=\mathcal{K}(\mathcal{P})(x) := \sum_{j=1}^{n} \left| \psi_j(x) \right|^2,\quad \forall x \in D.
    \label{christoffel}
\end{equation}
Now let $x_1, \ldots, x_m$ be drawn i.i.d.\ from a probability measure $\mu$, which is absolutely continuous with respect to $\rho$ and whose Radon--Nikodym derivative $d \mu / d \rho$ positive almost everywhere. Consider the $\tilde{f} \approx g$ approximation given by the weighted least-squares fit
\begin{equation}
\tilde{f} \in \operatorname{argmin}_{p \in \mathcal{P}} \frac{1}{m} \sum_{k=1}^m {w}(x_k)\, \bigl| f(x_k) - p(x_k) \bigr|^2 ,
\label{weighted-ls}
\end{equation}
where $w: D \rightarrow [0, \infty)$ is a weight function, given by $w = (d \mu / d \rho)^{-1}$.
In CS, one chooses $\mu$ as the \textit{CS measure}
\begin{equation}
    d\mu(x) = \frac{\mathcal{K}(x)}{n}\, d\rho(x)
    = \frac{\sum_{j=1}^n |\psi_j(x)|^2}{n}\, d\rho(x).
    \label{measurechris}
\end{equation}
This choice is near-optimal in the sense that if the number of samples satisfies $m \geq c\, n \log(n)$ for some positive constant $c$, the $\tilde{f}$ is a quasi-optimal approximation to $f$ from the subspace $\cP$.
For further details, we refer to \cite{adcock2025optimala}.

\section{Main algorithm: CAS-SRFE}
\label{sec:main}
We now develop CAS-SRFE.

\subsection{Christoffel Adaptive Sampling (CAS)}
\label{sec:CAS}
\begin{algorithm}[t]
\caption{Christoffel Adaptive Sampling (CAS)}
\label{alg:CAS}
\begin{algorithmic}
\STATE{\textbf{Input:}} Probability measure $\rho$ on a domain $D$, number of samples $m_i$, where $1 \leq m_1 < m_2 < \cdots$.
\STATE \textbf{Initialization:} Draw $x_1,...,x_{m_1} \sim_{\mathrm{i.i.d.}} \rho.$
\FOR{$i = 1, 2, \dots$}
  \STATE \textbf{Step 1:} Given data \( \left\{(x_k, f(x_k)\right\}_{k=1}^{m_i} \), compute an approximation \( \hat{f}_i \) and a subspace \( \mathcal{P}_i \) such that \( \hat{f}_i \in \mathcal{P}_i \). We require \( n_i=\dim(\mathcal{P}_i) < m_i<N \). 
  \STATE \textbf{Step 2:} Construct the {Christoffel function} \( \mathcal{K}: D \to \mathbb{R} \) of \( \mathcal{P}_i \),
   \[
    \mathcal{K}_i(x)= \mathcal{K}(\mathcal{P}_i)(x) := \sup \left\{\frac{|p_i(x)|^2}{\|p
    _i\|_{L_{\rho}^2(D)}}: p_i \in \mathcal{P}_i, p_i \neq 0\right\}, \quad \forall x \in D.
    \]
\STATE \textbf{Step 3:} Draw new samples \( x_{m_i+1}, \dots, x_{m_{i+1}} \sim_{\mathrm{i.i.d.}} \mu_i\), where $\mu_i$ is given by
  \[
    d\mu_i(x) = \frac{1}{n_i} \mathcal{K}_i(x) d\rho(x)
  \]
\ENDFOR
\end{algorithmic}
\end{algorithm}
As discussed, CS is a way to achieve near-optimal sampling for approximation in a linear space. CAS extends the applicability of CS to nonlinear approximation by adopting an adaptive strategy: it computes a sequence of approximations in adaptively-chosen subspaces, using CS at each step to select new samples. We now describe CAS in more detail.

CAS is summarized in Algorithm~\ref{alg:CAS}. It begins with $m_1$ MC samples, i.e. $x_1,\ldots,x_{m_1} \sim_{\mathrm{i.i.d.}} \rho$. These samples are then used to compute the initial approximation $\hat{f}_1$ via a (typically) nonlinear scheme, along with a subspace $\cP_1$ containing $\hat{f}_1$. The choice of $\cP_1$ is not fixed, and usually depends on the nonlinear scheme used to compute $\hat{f}_1$. See next for discussion. The key idea, however, is that we expect $\cP_1$ to be a good subspace of dimension $< m_1$ in which to approximate the target function $f$. Having done this, CAS then constructs the Christoffel function of this subspace along with the associated CS measure $\mu_1$. Finally, new samples samples $x_{m_1+1}, \ldots, x_{m_2}$ are drawn i.i.d.\ according to $\mu_1$. This process is then repeated iteratively.

As remarked, CAS is a general procedure for combining CS with adaptive sampling. In particular, the approximation $\hat{f}_i$ and the subspace $\cP_i$ in Step 1 can be constructed in different ways to yield different algorithms. In this work, we choose SRFE for Step 1. In particular, $\hat{f}_i$ is takes the form \eqref{RFM}, where at most $s = s_i$ of the coefficients are nonzero. We then use the features corresponding to these nonzero coefficients as a dictionary of functions to define the subspace $\cP_i$.

In the next subsections, we describe the specific implementations of Steps 1--3 of CAS in our CAS-SRFE framework.

\subsection{Implementing Step 1 of CAS via SRFE}
\label{sec:step1}
We use SRFE to compute $\hat{f}_i$. Fix $N \geq 1$, a feature function $\phi(x ; w)$, a probability distribution $\gamma$ for the features and draw $w_1,\ldots,w_N \sim_{\mathrm{i.i.d.}} \gamma$. Consider step $i$ of CAS and let
\begin{equation}
A^{(i)} = \frac{1}{\sqrt{m_i}} \left [ \phi(x_s ; w_r) \right ]_{s=1,\, r=1}^{m_i,\, N} \in \mathbb{C}^{m_i \times N},\quad b^{(i)}  = \frac{1}{\sqrt{m_i}} [f(x_s)]^{m_i}_{s=1} \in \mathbb{C}^{m_i}, 
\label{eq: Ax}
\end{equation} 
denote the random feature matrix and the vector of function samples, respectively, where $m_i$ is the number of samples at the $i$th iteration. Let $s_i < m_i$ be the sparsity level. SRFE involves finding a $s_i$-sparse approximate solution $\hat{c}^{(i)}$ to the linear system $A^{(i)} c \approx b^{(i)}$. In this work, we use OMP and HTP to construct $\hat{c}^{(i)}$.
These algorithms both produce $s$-sparse candidate solutions $\hat{c}$ to a linear system $A c \approx b$. They do so by iteratively computing $\hat{c}$ along with its support set $\hat{S}$. Here and throughout, we use the notation $\mathrm{supp}(c) \subseteq \{1,\ldots,N\}$ for the support set of a vector $c$, i.e., the set of indices corresponding to nonzero entries of $c$.

OMP, shown in Algorithm~\ref{alg:OMP}, finds an approximate sparse solution to by greedily adding indices to the support set. For a given $s$, OMP runs for $s$ iterations, adding a new index based on its correlation with the current residual.
HTP, shown in Algorithm~\ref{alg:HTP}, is a thresholding-based method that iteratively selects the $s$ largest components of $A^*(b - A c_{(k)})$ and solves a least-squares problem on the identified support ${S}$ until a suitable stopping criterion is met.

\begin{algorithm}[t]
\caption{Orthogonal Matching Pursuit (OMP)}
\label{alg:OMP}
\begin{algorithmic}

\STATE \textbf{Input:} Sparsity level $s$, matrix $A \in \mathbb{C}^{m\times N}$, vector $b \in \mathbb{C}^{m}$.
\STATE \textbf{Initialization:} Set $c_{(0)} = 0$, $S = \emptyset$.
\FOR{$k = 0, 1, \ldots, s-1$}
\STATE Select the index $j$ that maximizes the correlation with the current residual and update the support set:
\[
j \leftarrow \operatorname{argmax}_{j \in \{1,2,\ldots,N\}} \Big| \big(A^*(b - A c_{(k)}\big)_j \Big|, 
\quad S \leftarrow S \cup \{j\}.
\]
\STATE Solve the least-squares problem restricted to the selected support:
\[
c_{(k+1)} = \operatorname{argmin}_{z : \operatorname{supp}(z) \subseteq S} 
\| b - A z \|_2.
\]

\ENDFOR

\STATE \textbf{Output:} An $s$-sparse approximate solution $\hat{c} = c_{(s)}$ to the linear system $A c \approx b$ and its support set $\hat{S} = \mathrm{supp}(\hat{c})= S$. 
\end{algorithmic}
\end{algorithm}

\begin{algorithm}[htbp]
\caption{Hard Thresholding Pursuit (HTP)}
\label{alg:HTP}
\begin{algorithmic}

\STATE \textbf{Input:} Sparsity level $s$, matrix $A \in \mathbb{C}^{m\times N}$, vector $b \in \mathbb{C}^{m}$,
maximum iterations $\overline{k}$, tolerance $\texttt{tol}_{htp}$.
\STATE \textbf{Initialization:} Set $c_{(0)} = 0$, $S=\emptyset$.

\WHILE{$k = 0,1,2,\dots,\overline{k}-1$}
  \STATE Update the gradient and select the support:
\[ S \leftarrow L_s\big(c_{(k)} + A^*(b - A c_{(k)})\big),\]
  where $L_s$ selects the indices of the $s$ largest entries in magnitude.
  
  \STATE Solve the least-squares problem restricted to $S$:  
  \[c_{(k+1)} = \operatorname{argmin}_{z : \operatorname{supp}(z) \subseteq S} \| b - A z \|_2.\]
  
  \STATE Check stopping criterion:  
 $ \|c_{(k+1)} - c_{(k)}\|_2 < \texttt{tol}_{htp}.$
  If satisfied, terminate.
\ENDWHILE

\STATE \textbf{Output:} An $s$-sparse approximate solution $\hat{c}=c_{(k+1)}$ to the linear system $A c \approx b$ and its support set $\hat{S} = \mathrm{supp}(\hat{c})$.
\end{algorithmic}
\end{algorithm}

We make one modification to the standard OMP and HTP algorithms, which is column normalization. Column normalization often improves the performance of OMP and HTP by enforcing that the matrix used in either algorithm has to unit-norm columns. This is done by computing the Euclidean norm of each column and rescaling the matrix so that all columns have unit norm. Having computed $\hat{c}$, this process is then reversed to obtain a sparse approximate solution to the original linear system $A c \approx b$. This procedure is described in Algorithm~\ref{alg:colnorm}.  

\begin{algorithm}[t]
\caption{Column Normalization for OMP/HTP}
\label{alg:colnorm}
\begin{algorithmic}
\STATE \textbf{Input:} Matrix $A \in \mathbb{C}^{m \times N}$, vector $b \in \mathbb{C}^m$, sparsity level $s$ (with $\texttt{htp}_{\text{tol}}$ and $\Bar{k}$ for HTP).
\STATE For each column $j = 1, \dots, N$, compute the $\ell_2$ norm:
$\|A_j\| := \sqrt{\sum_{i=1}^m A_{ij}^2}.$
\STATE Form the normalized matrix:
$ A_{\text{norm}} \gets A \cdot \operatorname{diag}(1 / \|A_j\|)$
\STATE Apply OMP or HTP using $A_{\text{norm}}$, sparsity level $s$ and $b$ (with $\texttt{htp}_{\text{tol}}$ and $\Bar{k}$ for HTP) to obtain $\hat{c}$ and the support set $\hat{S}$.
\STATE Rescale the recovered coefficients:
$\hat{c} \gets \hat{c} \cdot \operatorname{diag}(1 / \|A_j\|)$
\STATE \textbf{Output:} Rescaled $s$-sparse approximate solution $\hat{c}$ to the linear system $A c \approx b$ and the support set $\hat{S} = \mathrm{supp}(\hat{c})$. 
\end{algorithmic}
\end{algorithm}

Having described OMP and HTP, we now return to Step 1 of Algorithm~\ref{alg:CAS}. When applied with $A = A^{(i)}$ and $b = b^{(i)}$ given by \eqref{eq: Ax} and some sparsity level $s = s_i$, both OMP and HTP output an $s_i$-sparse approximate solution $\hat{c}^{(i)}$ to $A^{(i)} c \approx b^{(i)}$ and its corresponding support set $\hat{S}^{(i)}$. Using this support set, we now conclude Step 1 by defining the subspace $\cP_i$ as
\begin{equation}
\mathcal{P}_i=\operatorname{span}\{\phi(\cdot,w_j)\mid j\in \hat{S}^{(i)}\}.
\label{eq:spansubspace}
\end{equation}
We then compute the Christoffel function of $\mathcal{P}_i$ in Step 2. We discuss this computation in the next subsection.

Before doing so, we now discuss one further modification to the above procedure, which is row reweighting. Row reweighting helps improve the numerical performance of CAS-SRFE by accounting for the change of measure that occurs when drawing samples from a probability measure that does not coincide with the underlying probability $\rho$. Let $x_1,\ldots,x_m \sim_{\mathrm{i.i.d.}} \mu$ and consider the random feature matrix $A = \frac{1}{\sqrt{m}} [ \phi(x_s ; w_r) ]^{m,N}_{s=1,r=1}$. Then
\[
\mathbb{E} ({A}^* {A})_{jk} = \mathbb{E} \left[ \frac{1}{m} \sum_{i=1}^{m} \phi(\cdot ; w_j) \phi(\cdot ; w_k)\right] = \int_{\mathbb{R}^d} \phi(x ; w_j) \phi(x ; w_k) d\mu(x).
\]
In reweighting, we scale the rows of $A$ so that this quantity coincides with the Gram matrix of the features $\{ \phi(\cdot ; w_j ) \}$. Specifically, let $w = (d \mu / d \rho)^{-1}$ and define
\[
\Bar{A} = \frac{1}{\sqrt{m}} \left [ \sqrt{w(x_s)} \phi(x_s ; w_r) \right ]^{m,N}_{s=1,r=1}.
\]
Then 
\begin{equation}
\mathbb{E} (\Bar{A}^* \Bar{A})_{jk} = \int_{\mathbb{R}^d} \phi(x ; w_j) \phi(x ; w_k) d\rho(x)
\label{exp-gram-prop}
\end{equation}
is now the Gram matrix of the features, as desired.

We implement row weighting in CAS-SRFE as follows. At the $i$th iteration, let $\mu_i$ denote the probability measure defined in Step 3 of Algorithm~\ref{alg:CAS} and let
\[
w_i : = \left (\frac{d \mu_i }{ d \rho} \right )^{-1} =  \left( \frac{1}{n_i} \mathcal{K}_i\right)^{-1}.
\]
Then we set
\begin{align*}
\Bar{A}^{(i)} &= \frac{1}{\sqrt{m_i}} \left [ \sqrt{w_i(x_s)} \phi(x_s, w_r) \right ]_{s,r=1}^{m_i,\, N} ,
\quad
 \Bar{b}^{(i)} &= \frac{1}{\sqrt{m_i}} \left [ \sqrt{w_i(x_s)} f(x_s) \right ]^{m_i}_{s=1},
\end{align*}
and use these in place of $A^{(i)}$ and $b^{(i)}$. In practice, since not all samples $\{x_k\}_{k=1}^{m_i}$ are drawn i.i.d.\ from $\mu_i$ in CAS-SRFE, the property \eqref{exp-gram-prop} may not exactly hold. Nonetheless, row weighting still generally improves its numerical performance.

\subsection{Implementing Step 2 of CAS}
\label{sec: step2}

Let $\cP_i$ be given by \eqref{eq:spansubspace}. To compute the Christoffel function of $\cP_i$ it suffices, in view of \eqref{christoffel}, to compute an orthonormal basis of $\cP_i$.  This can be done using the Gram matrix $G \in \mathbb{C}^{s_i\times s_i}$ of features, i.e.,
\begin{equation}
    G_{ij} = \int_D\overline{\phi(x; w_i)} \phi(x;w_j) \, d\rho(x),\quad i,j \in \hat{S}^{(i)}.
    \label{grammatrix}
\end{equation}
where $\overline{\phi(x;w_i)}$ is the complex conjugate of $\phi(x ;w_i)$.
Given the eigenvalues $\lambda_1 \geq \cdots \geq \lambda_{n_i} \geq 0$ and corresponding orthonormal eigenvectors $q_1, \dots, q_{n_i}$ of $G$, we obtain an orthonormal basis for $\mathcal{P}_i$ as follows:
\begin{equation}
\psi_j = \frac{1}{\sqrt{\lambda_j}} \sum_k (q_j)_k \phi(\cdot ; w_k).
\label{psi_j}
\end{equation}
Typically, the Gram matrix is ill-conditioned, since the features may have approximate linear dependencies, and so in practice, we discard all eigenvalues below a given threshold. Let $r_i \leq s_i$ be the number of eigenvalues above this threshold. This means that we compute the Christoffel function $\cK(\tilde{\cP_i})$, where 
\begin{equation}
\tilde{\cP_i}=\operatorname{span}\{\psi_1, \ldots, \psi_{r_i} \}, 
\label{tilde-Pi}
\end{equation}
rather than $\cK(\cP_i)$. The rationale behind this is that small eigenvalues correspond to nearly linearly dependent basis functions that contribute little to the approximation but can amplify numerical errors. Removing these basis functions improves the numerical stability of subsequent computations and downstream tasks. 

This aside, the key issue is computing the Gram matrix. In this paper, we consider two important scenarios where this is possible when using RFFs $\phi(x;w) = \exp(\mathrm{i} \langle x,w\rangle)$. The first is when $D=\mathbb{R}^d$ and $\rho=\mathcal{N}(0,\sigma^2 I)$. In this case,  we have
\begin{equation}
\label{gramplugin}
G_{ij}=\int_{\mathbb{R}^d}\phi^*(x;w_i)\phi(x;w_j)\,d\rho(x)
=\textrm{e}^{\mathrm{i}(w_j-w_i)^T\mu-\frac{1}{2}\sigma^2\|w_j-w_i\|^2_2}
\end{equation}
The second is the case where $D = [0,\infty)^d$ and $\rho = \mathrm{Exp}(\lambda)$ is the multivariate exponential distribution with parameter $\lambda > 0$. In this case,

$$
d \rho(x) = h(x) d x,\quad \text{where } h(x) = \begin{cases} \lambda \mathrm{e}^{-\lambda x} & x \geq 0 \\ 0 & x < 0 \end{cases}.
$$
In this case, Gram matrix has the form
\begin{equation}
G_{ij}=\int_{\mathbb{R}^d}\phi^*(x;w_i)\phi(x;w_j)\,d\rho(x)
=\prod_{k=1}^{d} \frac{\lambda}{\lambda - \mathrm{i} (w_j - w_i)_k}. 
\end{equation}
This matrix can also be computed explicitly in other cases, such as uniform random variables on cubes. In our experiments, we also consider cases where $\rho$ is a product of univariate Gaussian and exponential distributions with different parameters.

\subsection{Implementing Step 3 of CAS via MH}
\label{sec:MCMC}
\begin{algorithm}[t]
\caption{MH sampler for CAS Step 3}
\label{alg:MH}
\begin{algorithmic}

\STATE \textbf{Input:} Target density $p(x)$ (here $p(x)=p_i(x)=\dim(\mathcal{P}_i)^{-1}\,\mathcal{K}(\mathcal{P}_i)(x)\,h(x)$), initial point $x_0$, proposed $\sigma_1$ (the proposal density is $q(\cdot |x) = \mathcal{N}(\mu,\sigma_1 I_d)$), total kept samples $\texttt{total}$, burn-in $B$, thinning interval $T$.

\STATE \textbf{Initialization:} Set $x= x_0$, $S_p = \emptyset$, $t =0$, $\texttt{acc\_count} =0$.

\FOR{$t = 1,2,\dots$}
  \STATE  Compute the acceptance probability $\alpha = \min\!\left(1,\; \frac{p(x_{\text{prop}})}{p(x)}\right).$ Draw $u \sim \mathrm{Unif}(0,1)$ and set $x =  x_{\text{prop}}$ if $u \leq \alpha$ and $x = x$ otherwise.
  \STATE Burn-in and thinning:
  \IF{$t > B$ \textbf{ and } $(t-B) \bmod T = 0$}
     \STATE $S_p \gets S_p \cup \{x\}$; \quad $\texttt{acc\_count}  \gets \texttt{acc\_count} +1$.
     \IF{$k = \texttt{total}$}
        \STATE \textbf{break} \quad (Stop when the target number of samples is kept.)
     \ENDIF
  \ENDIF
\ENDFOR

\STATE \textbf{Post-processing:} Compute acceptance rate $\texttt{acc\_rate} \gets \texttt{acc\_count}/t$.

\STATE \textbf{Output:} Sample set $S_p=\{x_1,\dots,x_{\texttt{total}}\}$ drawn approximately from $p(x)$ and acceptance rate $\texttt{acc\_rate}$.
\end{algorithmic}
\end{algorithm}

We use an efficient MH sampler to implement Step 3 of CAS, a procedure that is described in Algorithm \ref{alg:MH}. In general,
MH provides an approach to generate samples that follow a given continuous multivariate distribution. Given $\mathcal{P}_i$ as defined in \eqref{eq:spansubspace}, the target density in Step 3 of CAS-SRFE is defined as 
\begin{equation}
p(x) =  p_i(x) : = \frac{1}{\dim(\mathcal{P}_i)} \mathcal{K}(\mathcal{P}_i)(x) h(x),
\label{targetdensity}
\end{equation}
where $\mathcal{K}$ as defined in \eqref{christoffel} and $h(x)$ is the density of $\rho$.
The algorithm generates a Markov Chain $x_1, x_2, \ldots$, where the new sample $x_{n+1}$ is only based on the previous state of the Markov chain, $x_{n}$. As more and more new samples are generated, the distribution of those samples will better approximate the target distribution $p(x)$. Therefore, the algorithm typically includes a burn-in phase to discard the initial samples, which may not yet follow the desired distribution closely \cite{robert2016metropolis}. 
In addition, the algorithm may include a {lag} parameter \(T\), also known as {thinning}. Thinning refers to a subsampling procedure in which only every \(T\)-th sample of the Markov chain is retained, i.e., $\{x_1, x_{1+T}, x_{1+2T}, \ldots\}$.
The purpose of thinning is to reduce the autocorrelation between successive samples, thereby improving their approximate independence, so that \(x_{n+T}\) is nearly independent of \(x_n\) \cite{robert1999monte}. Both burn-in and thinning are implemented in Algorithm \ref{alg:MH}.

MH will first generate a candidate $x^*$ from the proposal distribution $q(x^*|x_n)$, then calculate the acceptance probability $\alpha(x_n\rightarrow x^*)$:
 \[
    \alpha\left(x_n \rightarrow x^*\right) = \min \left(1, \frac{p\left(x^*\right)}{p\left(x_n\right)} \times \frac{q\left(x_n \mid x^*\right)}{q\left(x^* \mid x_n\right)}\right)
\]
In our case, we implement a Gaussian proposal distribution centered at the current state. This makes the proposal distribution symmetric: $q(x^* \mid x_n) = q(x_n \mid x^*)$.  As a result, the proposal ratio simplifies to
$
\frac{q(x_n \mid x^*)}{q(x^* \mid x_n)} = 1.
$
Therefore, the MH acceptance probability becomes:
\begin{equation}
\alpha(x_n \to x^*)
= \min\!\left(1,\,
\frac{p(x^*)}{p(x_n)} 
\cdot
\frac{q(x_n \mid x^*)}{q(x^* \mid x_n)}
\right)=\min\!\left(1,\,
\frac{p(x^*)}{p(x_n)}
\right).
\end{equation}
This simplified form is known as the Metropolis algorithm, which is a special case of the MH algorithm when the proposal distribution is symmetric.

Having computed the acceptance probability, the MH algorithm then generates a uniformly distributed random number $u$ over $[0,1]$ and determines the acceptance or rejection of the proposed sample by setting
    \[
    x_{n+1} =
    \begin{cases} 
    x^* & \text{if } u \leq \alpha\left(x_n \rightarrow x^*\right) \\ 
    x_n & \text{otherwise}.
    \end{cases}
    \]
In our work, the proposal distribution is Gaussian, $\mathcal{N}(x_n, \sigma_1 I_d)$, where the standard deviation $\sigma_1$ is adaptively tuned based on the acceptance rate. This process for tuning $\sigma_1$ shown in Algorithm~\ref{alg:MHtune}. The goal of this tuning algorithm is to find a value of variance $\sigma_1$ such that the acceptance rate falls within a relatively good range. This range is derived under strong theoretical assumptions in \cite{roberts1997weak}, but has been shown to be relatively robust in practice, even when the assumptions are not satisfied \cite{rosenthal2014optimising}. In this paper, we set the range to be $0.2-0.3$ for a faster tuning process.

We draw batches of samples and adjust $\sigma_1$ depending on whether the acceptance rate in each batch falls within a target range. If the acceptance rate exceeds the upper threshold, $\sigma_1$ is increased; otherwise, it is decreased. Although the tuning follows these rules, there is no rigorous theoretical justification for the relationship between the acceptance rate and the scale of $\sigma_1$. In addition, the bounds for the acceptance rate are determined empirically. Therefore, while tuning $\sigma_1$ is necessary, it does not require fine adjustment. For this reason, we run the tuning algorithm at the beginning of the main procedure to determine a uniform $\sigma_1$, which is then fixed throughout the subsequent sampling process.

\begin{algorithm}[t]
\caption{Adaptive Tuning of Proposal Scale $\sigma_1$}
\label{alg:MHtune}
\begin{algorithmic}
\STATE 
\textbf{Input:} Initial scale $\sigma_1 > 0$; target density $p(x)$; initial point $x_0$; burn-in $B$; thinning interval $T$.
\STATE \textbf{Initialization:} Set batch index $r = 1$, $\eta_{\mathrm{base}} = 0.05$, $\eta_{\mathrm{agg}} = 0.08$, $\texttt{batch\_size} = 200$, $k = B \bmod \texttt{batch\_size}$.
\STATE \textbf{Set acceptance bounds:}
\IF{$\mathrm{size}(x_0,2)=1$ (sample dimension $=1$), } 
    \STATE $\alpha_{\min} \gets 0.20$, $\alpha_{\max} \gets 0.26$
\ELSE
    \STATE $\alpha_{\min} \gets 0.35$, $\alpha_{\max} \gets 0.45$
\ENDIF
\STATE Calculate center $\alpha_{\mathrm{center}}=\frac{\alpha_{\min}+\alpha_{\max}}{2}.$
\FOR{$r = 1, 2, \dots,k$}
  \STATE \textbf{Step 1:} Using Algorithm~\ref{alg:MH} with current $\sigma_1$, target density $p(x)$, draw $\texttt{total}=\texttt{batch\_size}$ proposals (with burn-in $B$ and thinning $T$, initial point $x_0$) and output the acceptance rate for this batch $\texttt{acc\_rate}$.
  \STATE \textbf{Step 2:} Tune $\sigma_1$ based on $\texttt{acc\_rate}$:
  \IF{$\texttt{acc\_rate} > \alpha_{\max}$}
    \STATE $\Delta \gets \eta_{{\mathrm{base}}} + \eta_{{\mathrm{aggr}}}$
    \STATE $\sigma_1 \gets \sigma_1 \cdot (1 + \Delta)$
  \ELSIF{$\texttt{acc\_rate} < \alpha_{\min}$}
    \STATE $\Delta \gets \eta_{{\mathrm{base}}} + \eta_{{\mathrm{aggr}}}$
    \STATE $\sigma_1 \gets \sigma_1 \cdot (1 - \Delta)$
  \ELSE
    \IF{$\texttt{acc\_rate} > \alpha_{{\mathrm{center}}}$}
      \STATE $\sigma_1 \gets \sigma_1 \cdot (1 + 0.5\,\eta_{{\mathrm{base}}})$ 
    \ELSIF{$\texttt{acc\_rate} < \alpha_{{\mathrm{center}}}$}
      \STATE $\sigma_1 \gets \sigma_1 \cdot (1 - 0.5\,\eta_{{\mathrm{base}}})$ 
    \ENDIF
  \ENDIF
\ENDFOR
\STATE \textbf{Output:} The tuned proposal scale $\sigma_1$.
\end{algorithmic}
\end{algorithm}

\subsection{Boosting}
\label{sec:overview}

We have now described implementation of Steps 1--3 of CAS in the SRFE setting. In addition to these steps, we also implement a boosting procedure based on \cite{haberstich2022boosted} to enhance the performance of the algorithm. 

Boosting is a way to improve standard CS for linear least-squares approximation. It does so by repeatedly drawing sets of sample points and then selecting the best, according to an empirical criterion. To explicate, consider the setup of \S \ref{sec:cs}. Let $\cP \subseteq L^2_{\rho}(D)$ be an $n$-dimensional subspace, $\{ \psi_j \}^{n}_{j=1}$ be an orthonormal basis for $\cP$, $w(x)$ be a weight function and $x_1,\ldots,x_m \in D$ be a set of sample points. Note that, for the moment, we do not assume the $x_i$ are drawn from a probability measure $\mu$, nor that $w$ is the reciprocal of the Radon-Nikodym derivative of $\mu$.  As discussed in, e.g., \cite{adcock2025optimala} the accuracy and stability of the weighted least-squares fit \eqref{weighted-ls} are determined by the \textit{discrete stability constant}, which is the minimal singular value $\sigma_{\min}(B)$ of the least-squares matrix
\[
B = \left [ \sqrt{w(x_i)} \psi_j(x_i) \right ]^{m,n}_{i,j=1} \in \mathbb{C}^{m \times n},
\]
with a larger $\sigma_{\min}(B)$ implying better accuracy and stability.
Hence, in boosting, we repeatedly generate candidate sets of sample points, then pick the one that maximizes $\sigma_{\min}(B)$. For more details on boosting in the context of CS, see \cite{haberstich2022boosted}.

The boosting procedure for CAS is described in Algorithm~\ref{alg:boosted_single}. In it, we initialize the stability constant $\alpha_{\mathrm{best}}$ to $0$ and run the procedure for a prescribed number of iterations. At iteration $i$, we first draw $m_{i+1} - m_i$ new samples from the measure $\mu_{i}$ and combine them with the previously drawn samples $\{x_k\}_{k=1}^{m_i}$ to obtain the enlarged set $X = \{x_k\}_{k=1}^{m_{i+1}}$. From this updated sample set, we form weighted the least-squares matrix $B$. Note that the orthonormal basis needed for this step has already been computed in Step 2 of CAS-SRFE in order to compute the Christoffel function. Next, we compute the smallest singular value of $B$ and compare it with the current value of $\alpha_{\mathrm{best}}$; if it is larger, we update $\alpha_{\mathrm{best}}$ accordingly. We then repeat this procedure until the prescribed number of boosting steps is reached.

\begin{algorithm}[t]
\caption{{Boosting in CAS}}
\label{alg:boosted_single}
\begin{algorithmic}
\STATE{\textbf{Input:}} Number of boosting iterations $Bo$; prior sample set $X=\{x_k\}^{m_{i}}_{k=1}$; target distribution $p(x)$ (sampling measure $\mu_i$), initial guess for proposed scale $\sigma_1$, initial point $x_0$, burn-in $B$ and thinning interval $T$; orthonormal basis $\{\psi_j\}_{j=1}^{r_i}$ for $\tilde{\cP}_i$, where $\tilde{\cP}_i$ is as in \eqref{tilde-Pi}; 

\STATE{\textbf{Initialization:}} Set best stability value $\alpha_{\text{best}} = 0$.

\FOR{$t = 1, 2, \dots, Bo$}
\STATE Draw new samples \( x_{m_{i}+1},\ldots,x_{m_{i+1}} \sim_{\mathrm{i.i.d.}} \mu_i \) using Algorithm~\ref{alg:MH} with $\sigma_1$ selected via Algorithm~\ref{alg:MHtune}, with target density $p(x)$, initial point $x_0$, burn-in $B$ and thinning $T$ for both algorithms (initial guess for $\sigma_1$ for Algorithm~\ref{alg:MHtune}).

  \STATE Compute the weighted least-squares matrix
\[
B =  \frac{1}{\sqrt{m_i}} \left [ \sqrt{w_i(x_s)}\psi_r(x_s) \right ]^{m_{i+1},r_i}_{s,r=1} \in \mathbb{C}^{m_{i+1} \times r_i},\quad \text{where }w_i = (d \mu_i / d \rho)^{-1}.
\]
  \STATE Compute the discrete stability constant $\alpha = \sigma_{\min}(B)$
  \STATE If \( \alpha > \alpha_{\text{best}} \), update:
  $\alpha_{\text{best}} \gets \alpha, \quad X_{\text{best}} \gets \left\{x_k\right\}_{k=1}^{m_{i+1}}$
\ENDFOR
\STATE{\textbf{Output:}} Selected sample set \( X_{\text{best}} \) with maximal discrete stability constant \( \alpha_{\text{best}} \).
\end{algorithmic}
\end{algorithm}

\subsection{Overview of CAS-SRFE}
This completes the description of the key steps of CAS in the case of SRFE. We now summarize the CAS-SRFE framework in Algorithm \ref{alg:CAS-SRFE}, including all supporting algorithms.

\begin{algorithm}[t]
\caption{CAS-SRFE}
\label{alg:CAS-SRFE}
\begin{algorithmic}
\STATE \textbf{Input:} 
Number of samples $1 \leq m_1 < m_2 < \dots$; underlying measure $\rho$; number of boosting steps $Bo$; eigenvalue truncation tolerance \texttt{tol}; number of random features $N$; measure $\gamma$; HTP parameters: maximum iterations $\overline{k}$, tolerance $\texttt{tol}_{htp}$;
MH parameters: initial point $x_0 \in D$, initial guess for the proposal scale $\sigma_1$, thinning interval $T$, burn-in $B$.

\STATE \textbf{Initialization:} 
Set $S = \{1,2,\ldots,\lfloor m_1/4 \rfloor\}$.

 \STATE \textbf{Precompute:}
Generate frequencies $w_1, \ldots, w_N \sim_{\mathrm{i.i.d.}} \gamma$ and define the corresponding  
random features $\phi(\cdot; w_j)$, $j = 1,\ldots,N$.

\STATE \textbf{Tuning of $\sigma_1$:}
Using the initialized support set $S$, calculate the target density $p(x)$. Run Algorithm~\ref{alg:MHtune} with target density $p(x)$, initial proposed scale $\sigma_1$, burn-in $B$, and thinning interval $T$ to obtain the tuned value of $\sigma_1$ for subsequent MH sampling. 
\FOR{$i = 1,2,\ldots$}

\IF{$i=1$}
\STATE Draw $x_1, \ldots, x_{m_1} \sim_{\mathrm{i.i.d.}} \rho$.
\ELSE
 \STATE \textbf{Step 1a: Run sparse recovery algorithm for SRFE.} Construct $A = A^{(i)}$ and $b = b^{(i)}$ as in \eqref{eq: Ax} and set $s_i = \lfloor m_i/4 \rfloor$. Run Algorithm~\ref{alg:colnorm}, using either OMP  (Algorithm~\ref{alg:OMP})  or HTP (Algorithm~ \ref{alg:HTP}) with $s$, $A$ and $b$ (also $\texttt{tol}_{htp}$ and $\overline{k}$ if using HTP). Obtain the $s$-sparse approximate solution $\hat{c}_i$ and support set $\hat{S}_i$.

   \STATE \textbf{Step 1b: Subspace construction.}
    Construct $\hat f^{(i)}$ via~\eqref{RFM} with $\hat{c} = \hat{c}^{(i)}$ and define
    $
    \mathcal P_i = \operatorname{span}\{\phi(\cdot,w_j): j \in \hat{S}_i\}.
    $

\STATE \textbf{Step 2a: Basis construction.} 
Compute an orthonormal basis $\{\psi_k\}_{k=1}^{r_i}$ of $\tilde{\cP_i} =\operatorname{span}\{\psi_1, \ldots, \psi_{r_i} \} \subseteq \cP_i$ by computing eigendecomposition of the Gram matrix~\eqref{grammatrix} and using \eqref{psi_j}. Here $r_i$ the number of eigenvalues exceeding \texttt{tol}.

\STATE \textbf{Step 2b: Compute Christoffel function.} Compute $\mathcal K(\tilde{\mathcal P}_i)$ via~\eqref{christoffel} with $n = r_i$ and define the sampling measure $\mu_i$ via~\eqref{measurechris}.

\STATE \textbf{Step 3: Sampling via MH with boosting.} Run Algorithm~\ref{alg:boosted_single} with $Bo$, $X = \{x_k\}_{k=1}^{m_i}$, target density $p(x)$ (sampling measure $\mu_i$), all MH parameters, $\sigma_1$ and orthonormal bais $\{\psi_k\}_{k=1}^{r_i}$ to obtain the new sample set $\{ x_k \}^{m_{i+1}}_{k=1}$.
\ENDIF
\ENDFOR
\end{algorithmic}
\end{algorithm}

\section{Numerical Experiments}
\label{sec:numerics}

We now present a series of numerical experiments to demonstrate the performance of CAS-SRFE. 

\subsection{General setup} 
We first describe the general setup for our experiments, including parameter choices. For the main algorithm, we use $Bo=5$ boosting iterations and a truncation tolerance of $\texttt{tol}=1e-10.$ For the MH algorithm, we use a thinning interval of $T=15$, a burn-in of $B=5000$ and the initial point $x_0=0$ unless otherwise specified. For tuning the proposal scale, we use $\eta_{\mathrm{base}}=0.05$, $\eta_{\mathrm{aggr}}=0.08$ and a batch size of $M_b=200$. If $d=1$, then we use $[\alpha_{\min}, \alpha_{\max}]=[0.35, 0.45]$, while if $d>1$ we use $[\alpha_{\min}, \alpha_{\max}]=[0.2, 0.26]$. In both cases, we set $\alpha_{\mathrm{center}}=(\alpha_{\min}+\alpha_{\max})/2$. For HTP, we set $\texttt{htp}_{\text{tol}}=1e-12$ and $\Bar{k}=100$.

For the SRFE, we use RFFs with $\phi(x;w) = \exp(\mathrm{i} \langle x , w \rangle)$. We use $\gamma = \mathcal{N}(0, \sigma^2_w I)$ to generate the random features. Unless otherwise specified, we fix $\sigma_w=1$. We use several different choices for the underlying distribution $\rho$, specified later in the experiments. We compare the CAS-SRFE approximation against the standard SRFE approximation computed from nonadaptive MC samples drawn i.i.d.\ from $\rho$. We refer to the former as CAS and the latter as Non-Adaptive Sampling (NAS) in all figures.
For testing, we randomly draw $K_{t}=10,000$ points according to $\rho$. If $\hat{f}$ denotes an approximation to $\hat{f}$, we use these points to approximately evaluate the relative $L^2_{\rho}$-norm error $\|f-\hat{f}\|_{L^2_{\rho}}/\|f\|_{L^2_{\rho}}$.
We run multiple trials and, following \cite[\S A.1]{adcock2022sparse} report the (geometric) mean and standard deviation. For all of our experiments, we fix the number of trials $\texttt{trials}=30$. For the parametric differential equations in \S \ref{sec:parametric}, we solve the systems using $\texttt{ode45}$ for testing the approximation results.

\subsection{MH algorithm diagnostics}

In this section, we present numerical results that demonstrate the behaviour and diagnostics of the MH sampling algorithm. We run experiments with the underlying distribution being the normal distribution and the exponential distribution, illustrating the agreement between the drawn samples and the target density, as well as how the Christoffel function modifies the underlying distribution. For the one-dimensional cases, we also present trace plots to evaluate the mixing of the Markov chain.

\subsubsection*{Normal distribution examples} 
\begin{figure}[htbp]
    \centering
    \begin{subfigure}[t]{0.48\textwidth}
        \centering
       \includegraphics[
            width=\linewidth,
            trim={0mm 0mm 0mm 0mm},
            clip
        ]{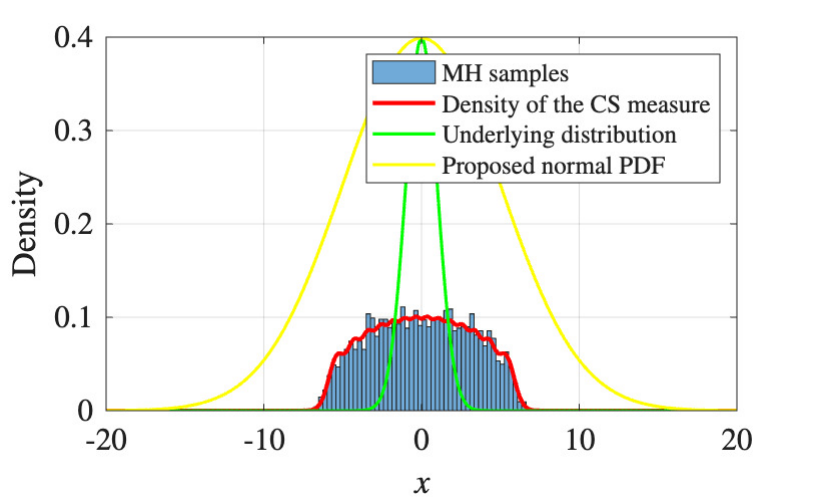}
        \caption{MH sampling for the 1D CS measure (histogram vs. target density). The reference density and proposed density are included for context.
        }
        \label{fig:1dmh_samp}
    \end{subfigure}
    \hfill
    \begin{subfigure}[t]{0.48\textwidth}
        \centering
         \includegraphics[
            width=\linewidth,
            trim={0mm 0mm 0mm 0mm},
            clip
        ]{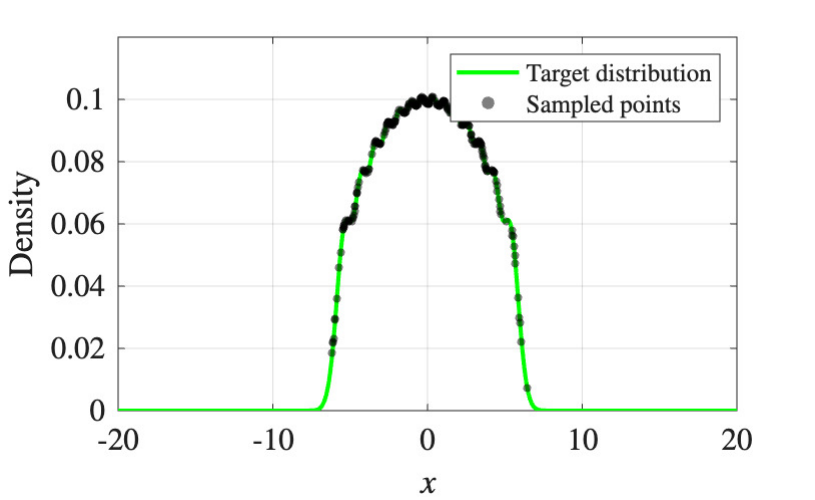}
        \caption{CS target density (green) with MH samples (black dots) overlaid.}
        \label{fig:1d_christoffel}
    \end{subfigure}
    \begin{subfigure}[t]{0.7\textwidth}
        \centering
         \includegraphics[
            width=\linewidth,
            trim={0mm 0mm 0mm 0mm},
            clip
        ]{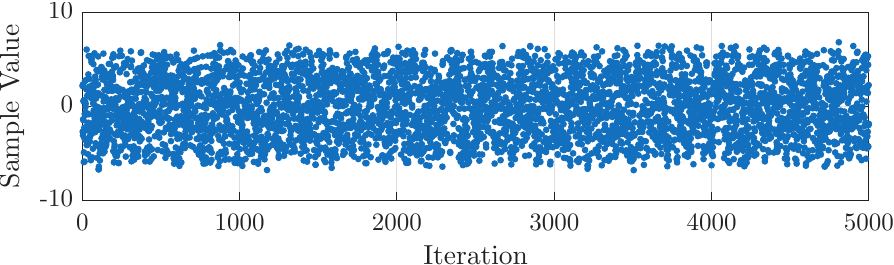}
        \caption{Trace plot of the MH chain (1-D), i.e., sample values over iterations.}
        \label{fig:1d_trace}
    \end{subfigure}
    \caption{Diagnostics of the MH sampling algorithm with {the CS measure \eqref{measurechris}} with $\rho = \mathcal{N}(0,\sigma^2 I)$. The parameters are $s=100$ (number of random features used for computing the Christoffel function), $d = 1$, $\sigma=1$, $\sigma_1=5$, $B=2000$, total samples \texttt{total} $=5000$, thinning $T=5$. 
    }
    \label{fig:1dmh}

\end{figure}

\begin{figure}
    \centering
    \begin{subfigure}[t]{0.45\textwidth}
        \centering
        \includegraphics[
            width=\linewidth,
            clip
        ]{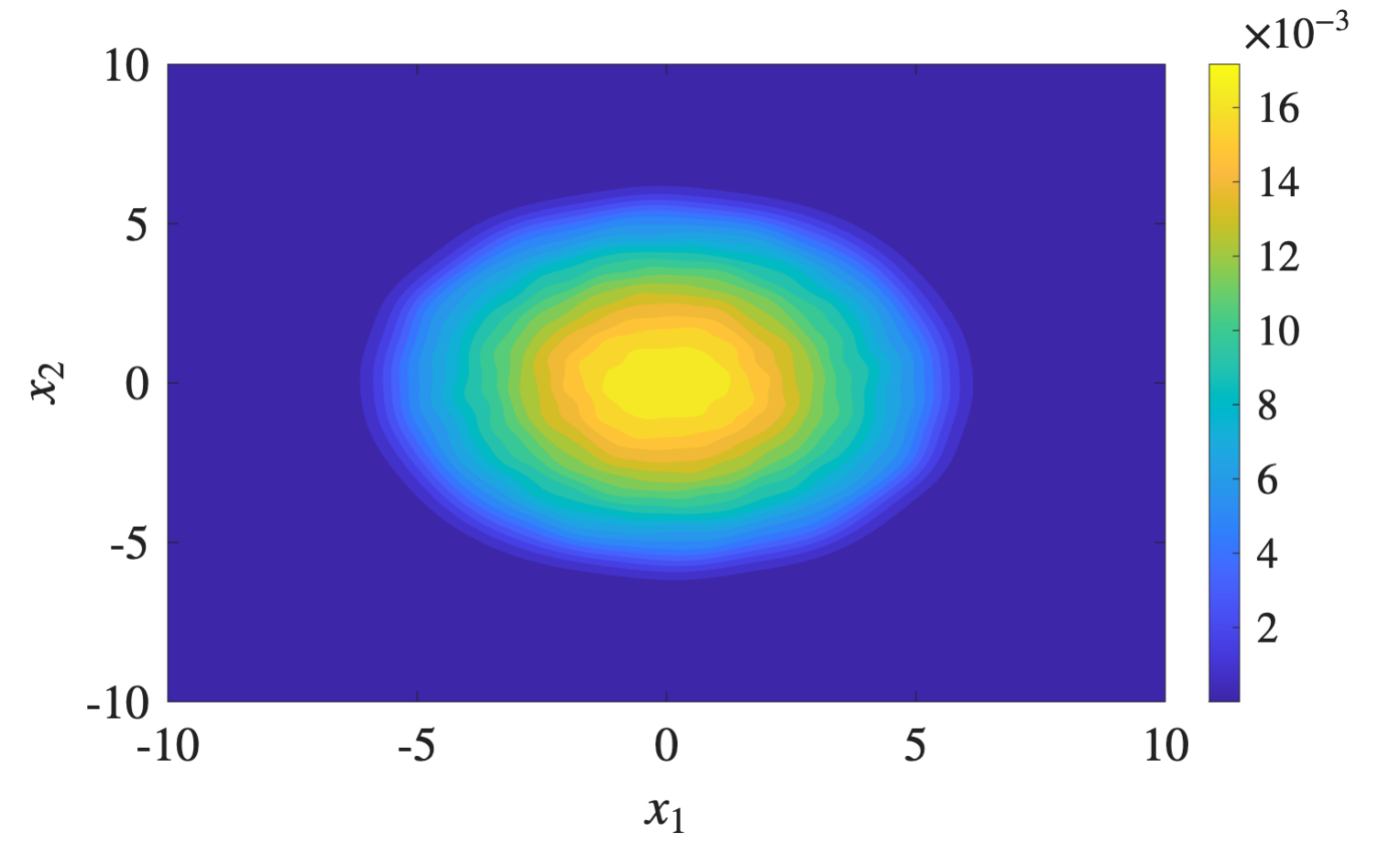}
        \caption{Density of the CS measure.}
        \label{fig:mh2d_true_bx}
    \end{subfigure}
    \hfill
    \begin{subfigure}[t]{0.45\textwidth}
        \centering
        \includegraphics[
            width=\linewidth,
            clip
        ]{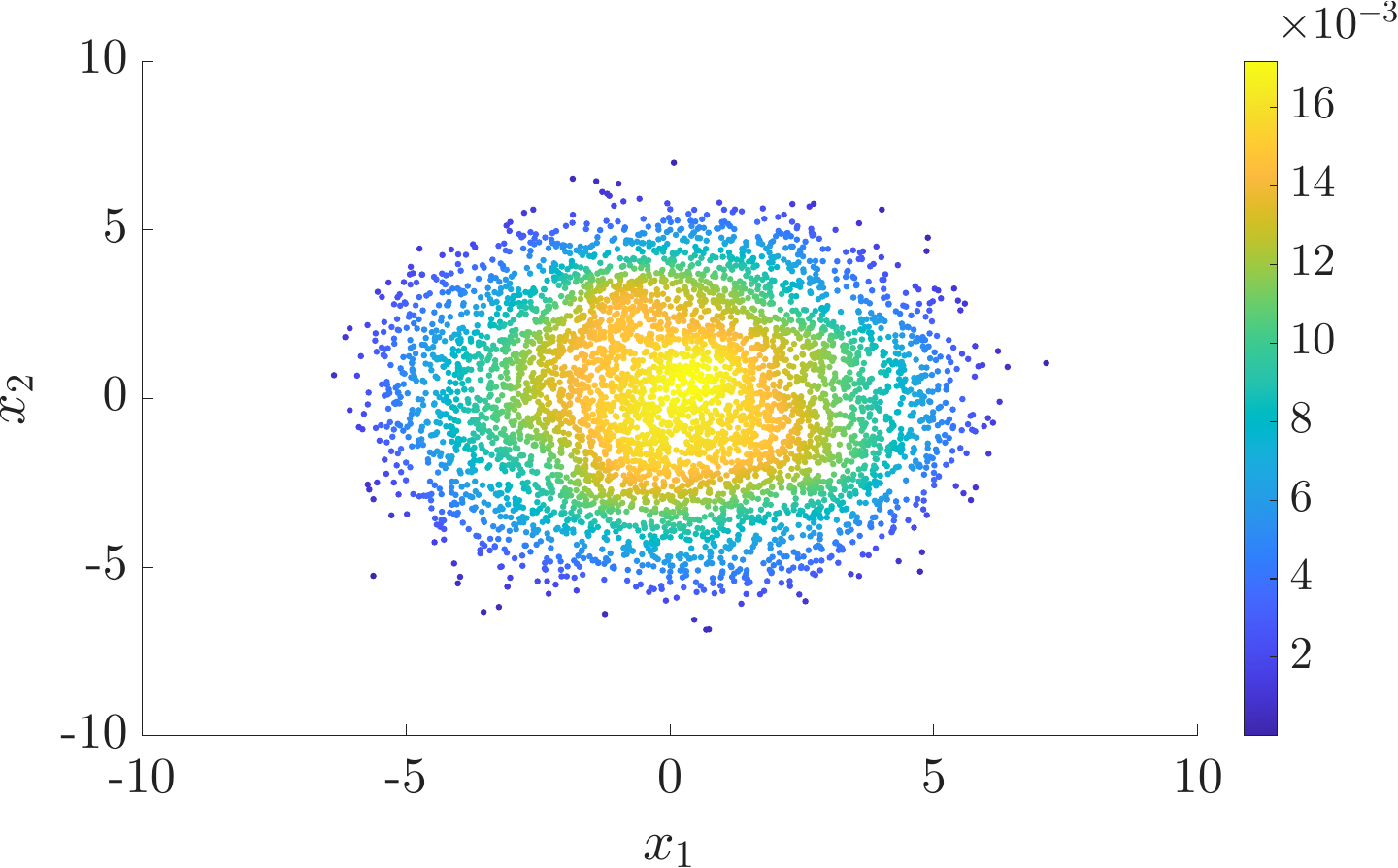}
        \caption{MH samples coloured by density.}
        \label{fig:mh2d_samples_colored}
    \end{subfigure}
    \begin{subfigure}[t]{0.45\textwidth}
        \centering
        \includegraphics[
            width=\linewidth,
            clip
        ]{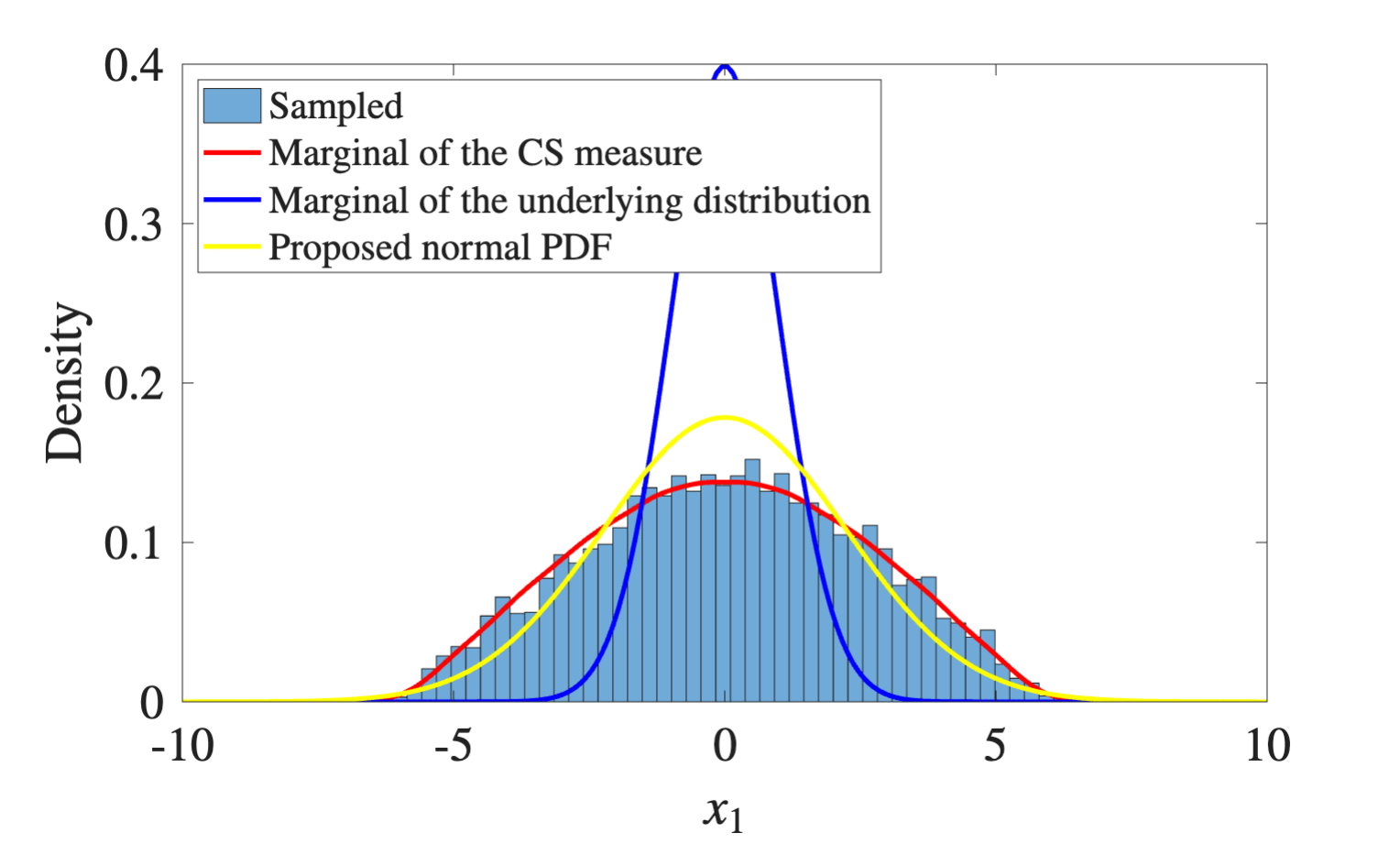}
        \caption{Marginal distribution of $x_1$.}
        \label{fig:mh2d_marginal_x1}
    \end{subfigure}
    \hfill
    \begin{subfigure}[t]{0.45\textwidth}
        \centering
        \includegraphics[
            width=\linewidth,
            clip
        ]{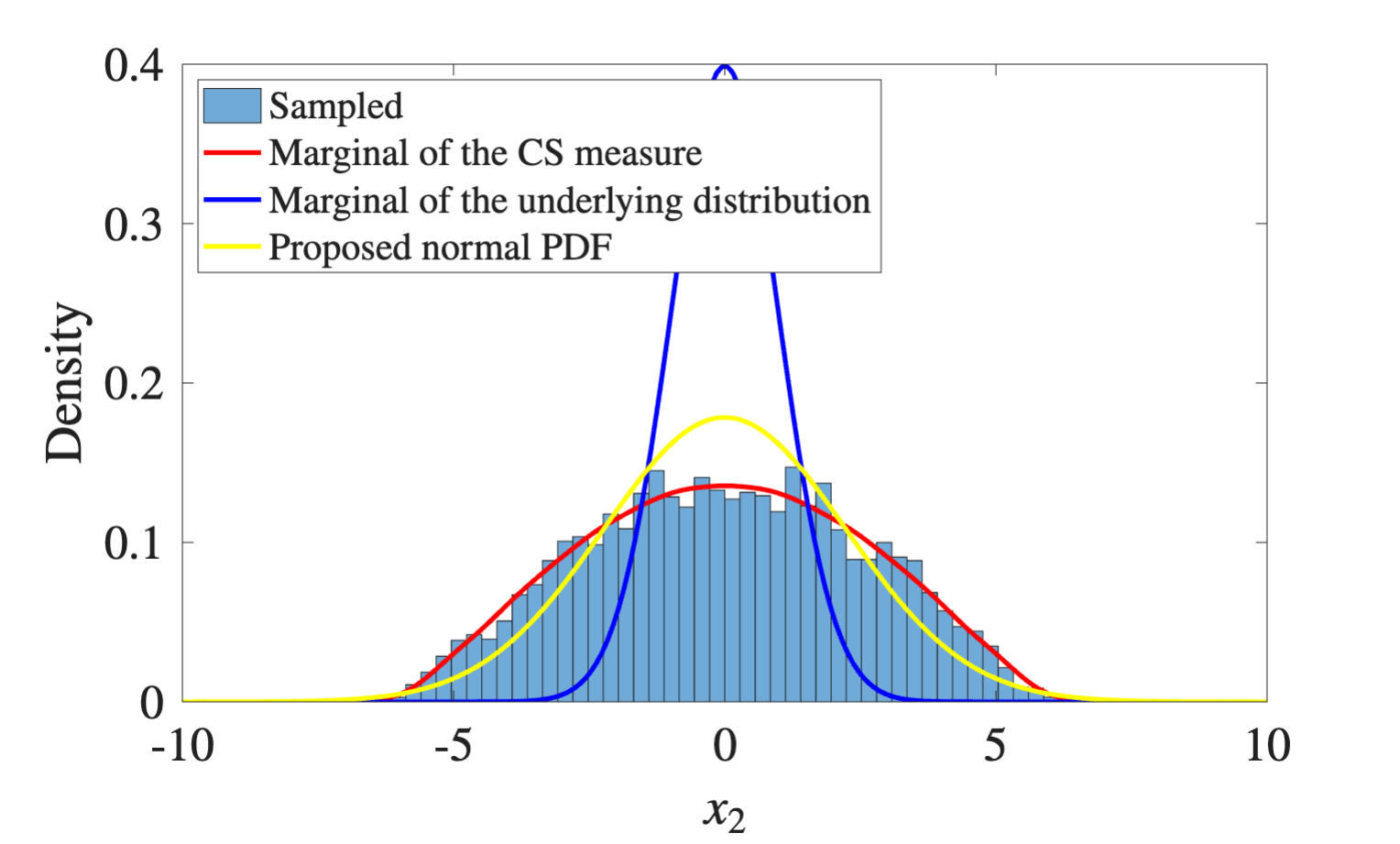}
        \caption{Marginal distribution of $x_2$.}
        \label{fig:mh2d_marginal_x2}
    \end{subfigure}
    \caption{Diagnostics of the MH sampling algorithm with the CS measure \eqref{measurechris} with $\rho = \mathcal{N}(0,\sigma^2 I)$. The parameters are s = 100, $d$ = 2, $\sigma$ = 1, $\sigma_1$ = 5, burn-in $B = 5000$, total samples \texttt{total} = 5000, thinning  $T= 15$.}
\label{fig:2dmh}
\end{figure}

Figures~\ref{fig:1dmh} and \ref{fig:2dmh} show simulations of the MH sampling algorithm when the underlying distribution $\rho = \mathcal{N}(0,\sigma^2 I)$. Figures~\ref{fig:1dmh}(a),(b) plot the density corresponding to the CS measure \eqref{measurechris} based on a subspace $\cP$ spanned by $s = 100$ random features, and compare it with the empirical distribution of the drawn samples, showing a close match. Figure~\ref{fig:1dmh}(c) displays a trace plot, which evaluates the mixing of the chain. The mixing time of a Markov chain is the time required for the chain to approach its steady-state distribution. An ideal trace plot should converge to a stationary region and fluctuate around the target density without noticeable trends, as observed here.
Figure~\ref{fig:2dmh} illustrates the two-dimensional case. Once again, the patterns match well, indicating that the samples are following the target density properly. The bottom figures show a comparison of the marginal distributions in two separate dimensions with the true distribution.


\subsubsection*{Exponential distribution examples}
\begin{figure}[htbp]
    \centering
    \begin{subfigure}[t]{0.48\textwidth}
        \centering
        \includegraphics[
            width=\linewidth,
            trim={0mm 0mm 0mm 0mm},
            clip
        ]{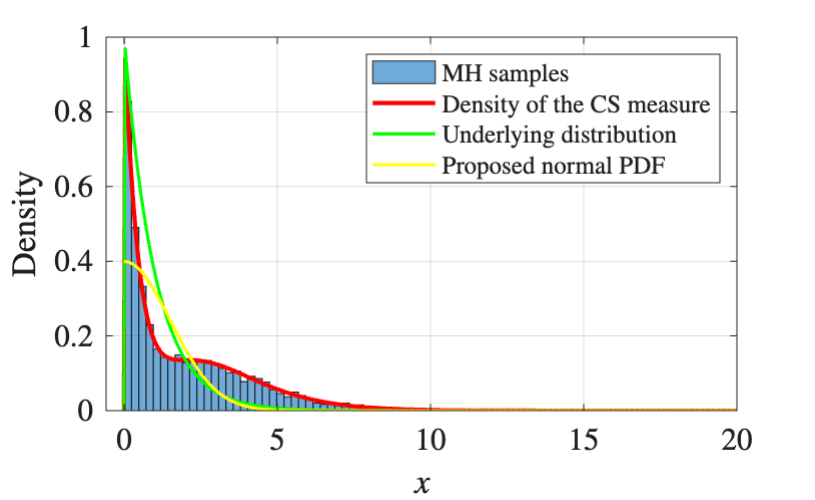}
           \caption{MH sampling for the 1D CS measure (histogram vs. target density). The reference density and proposed density are included for context. }
        \label{fig:1dmh_samp_exp}
    \end{subfigure}
    \hfill
    \begin{subfigure}[t]{0.48\textwidth}
        \centering
        \includegraphics[
            width=\linewidth,
            trim={0mm 0mm 0mm 0mm},
            clip
        ]{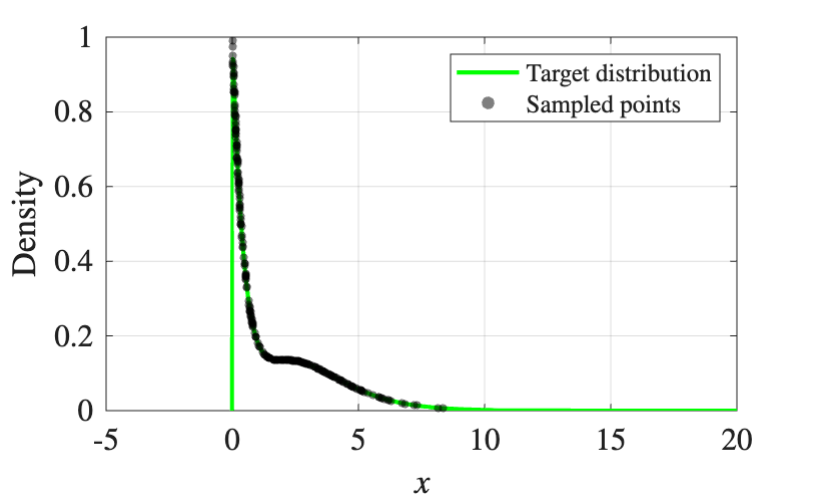}
           \caption{CS target density (green) with MH samples (black dots) overlaid.}
        \label{fig:1d_christoffel_exp}
    \end{subfigure}
    \vspace{0.6em}
    \begin{subfigure}[t]{0.7\textwidth}
        \centering
        \includegraphics[
            width=\linewidth,
            trim={0mm 0mm 0mm 0mm},
            clip
        ]{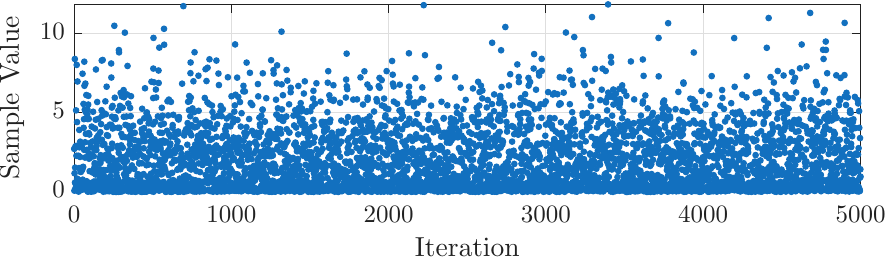}
          \caption{Trace plot of the MH chain (1-D), i.e., sample values over iterations.}
        \label{fig:1d_trace_exp}
    \end{subfigure}  
    \caption{Diagnostics of the MH sampling algorithm with the CS measure \eqref{measurechris} with $\rho = \mathrm{Exp}(\lambda)$, $\lambda=1e-3$. The parameters are $s = 100$, $d$ = 1, $\sigma$ = 1, $\sigma_1$ = 1.5, $\sigma_w=10^{-3}$, burn-in $B = 5000$, total samples $\texttt{total}= 5000$, thinning $T= 30$.}
    \label{fig:1dmh_exp}
\end{figure}

\begin{figure}[htbp]
    \centering
    \begin{subfigure}[t]{0.45\textwidth}
        \centering
        \includegraphics[
            width=\linewidth
        ]{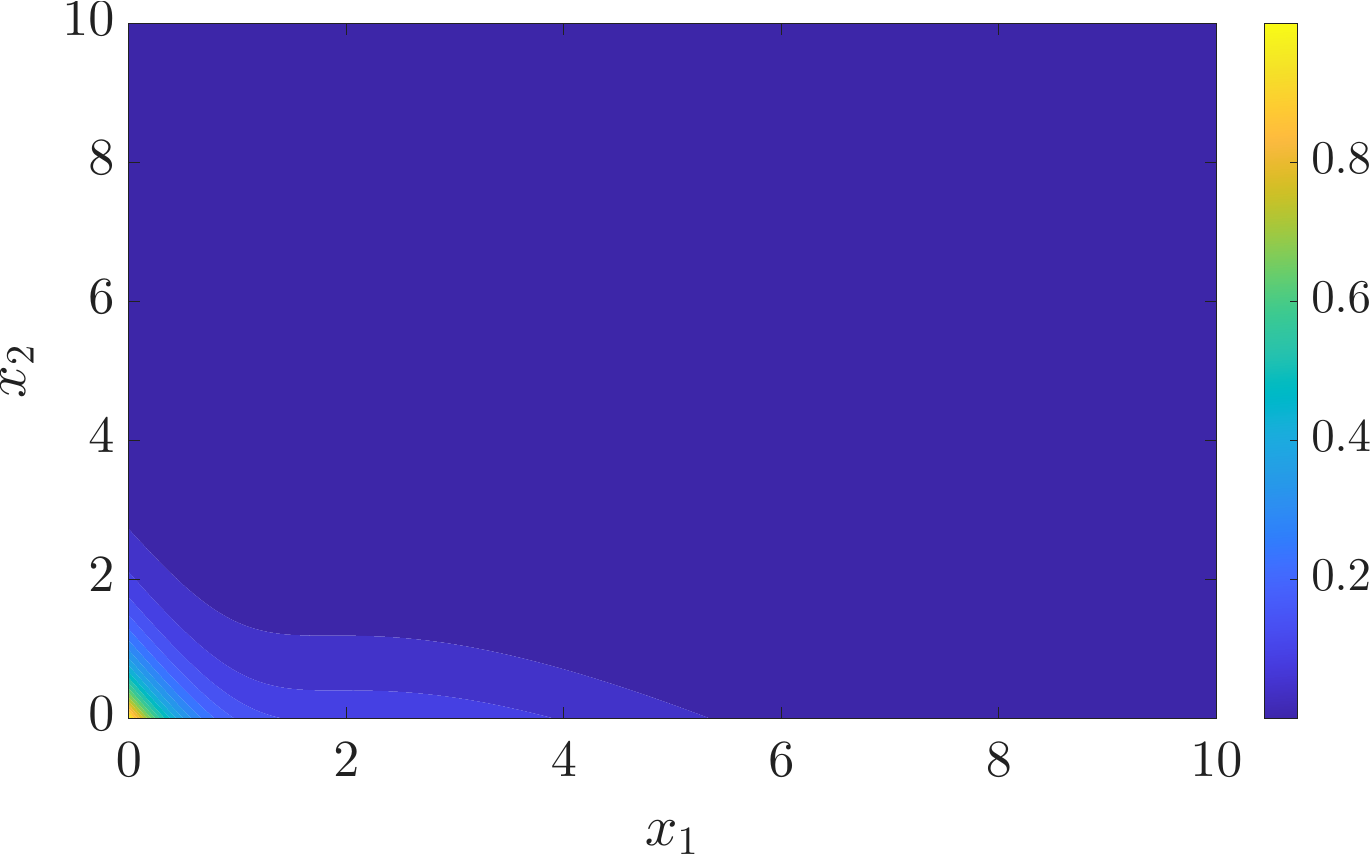}
        \caption{Density of the CS measure.}
        \label{fig:mh2d_true_bx_exp}
    \end{subfigure}
    \hfill
    \begin{subfigure}[t]{0.45\textwidth}
        \centering
        \includegraphics[
            width=\linewidth
        ]{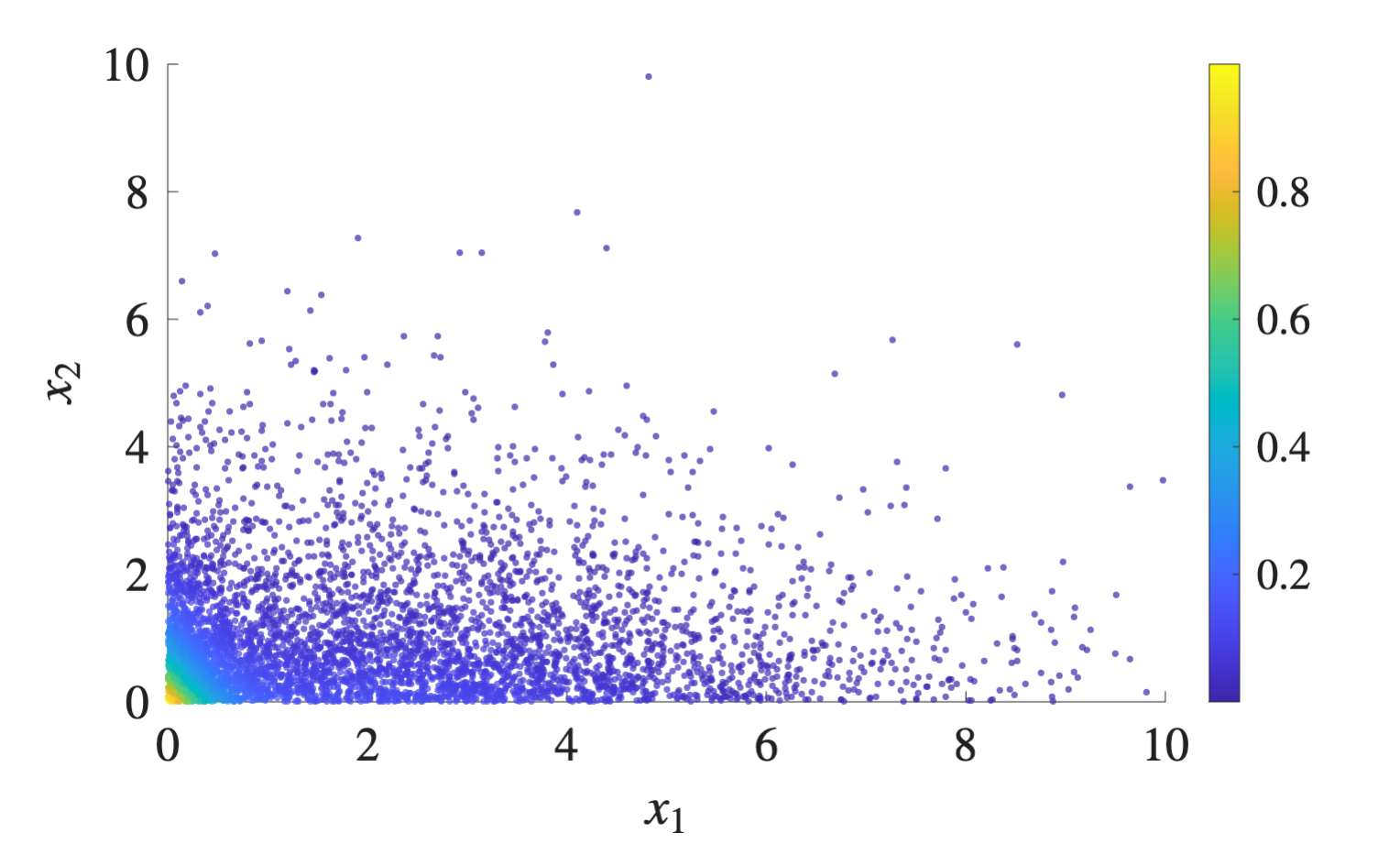}
        \caption{MH samples coloured by density.}
        \label{fig:mh2d_samples_colored_exp}
    \end{subfigure}
    \begin{subfigure}[t]{0.45\textwidth}
        \centering
        \includegraphics[
            width=\linewidth,
            clip
        ]{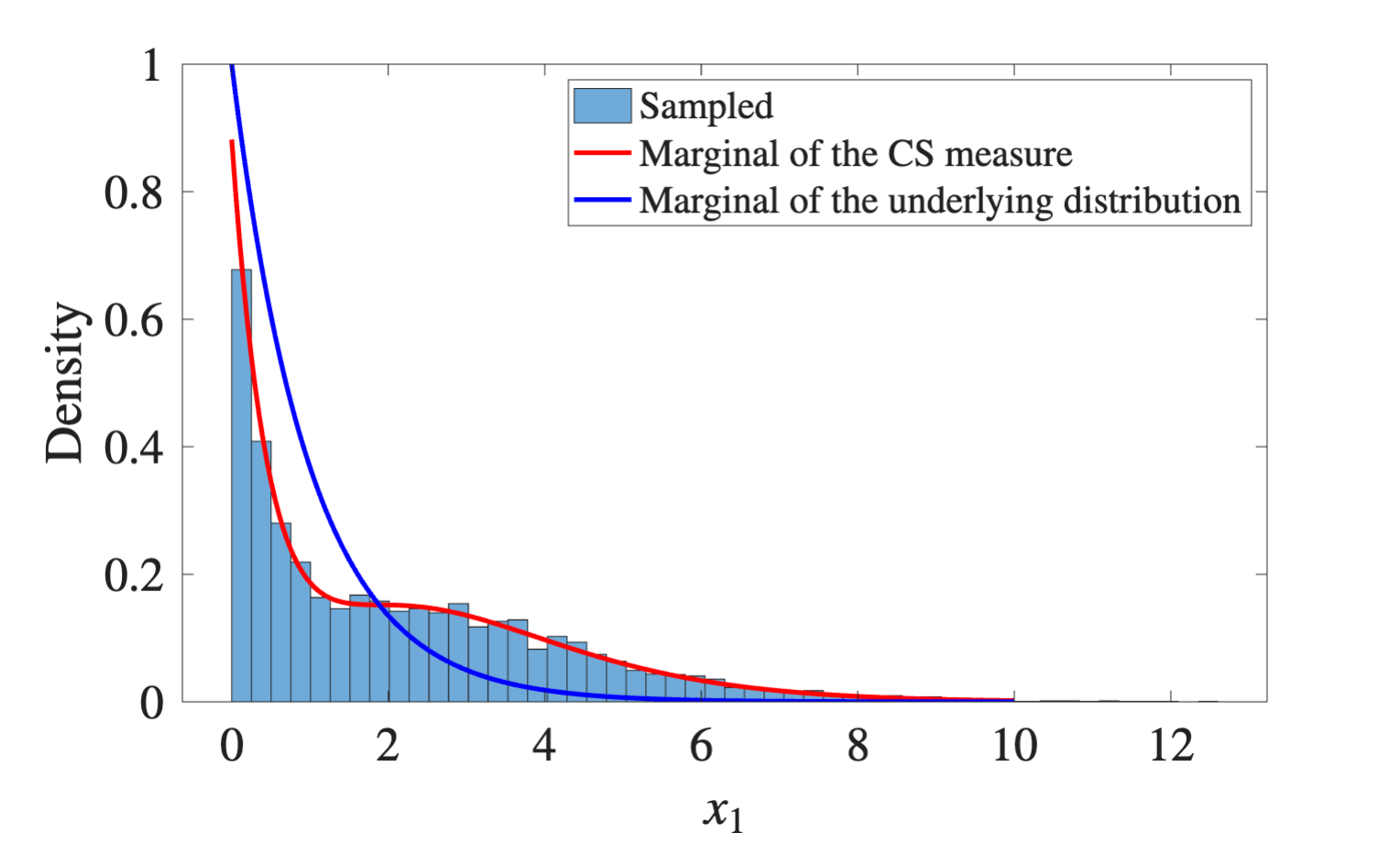}
        \caption{Marginal distribution of $x_1$.}
        \label{fig:mh2d_marginal_x1_exp}
    \end{subfigure}
    \hfill
    \begin{subfigure}[t]{0.45\textwidth}
        \centering
        \includegraphics[
            width=\linewidth,
            clip
        ]{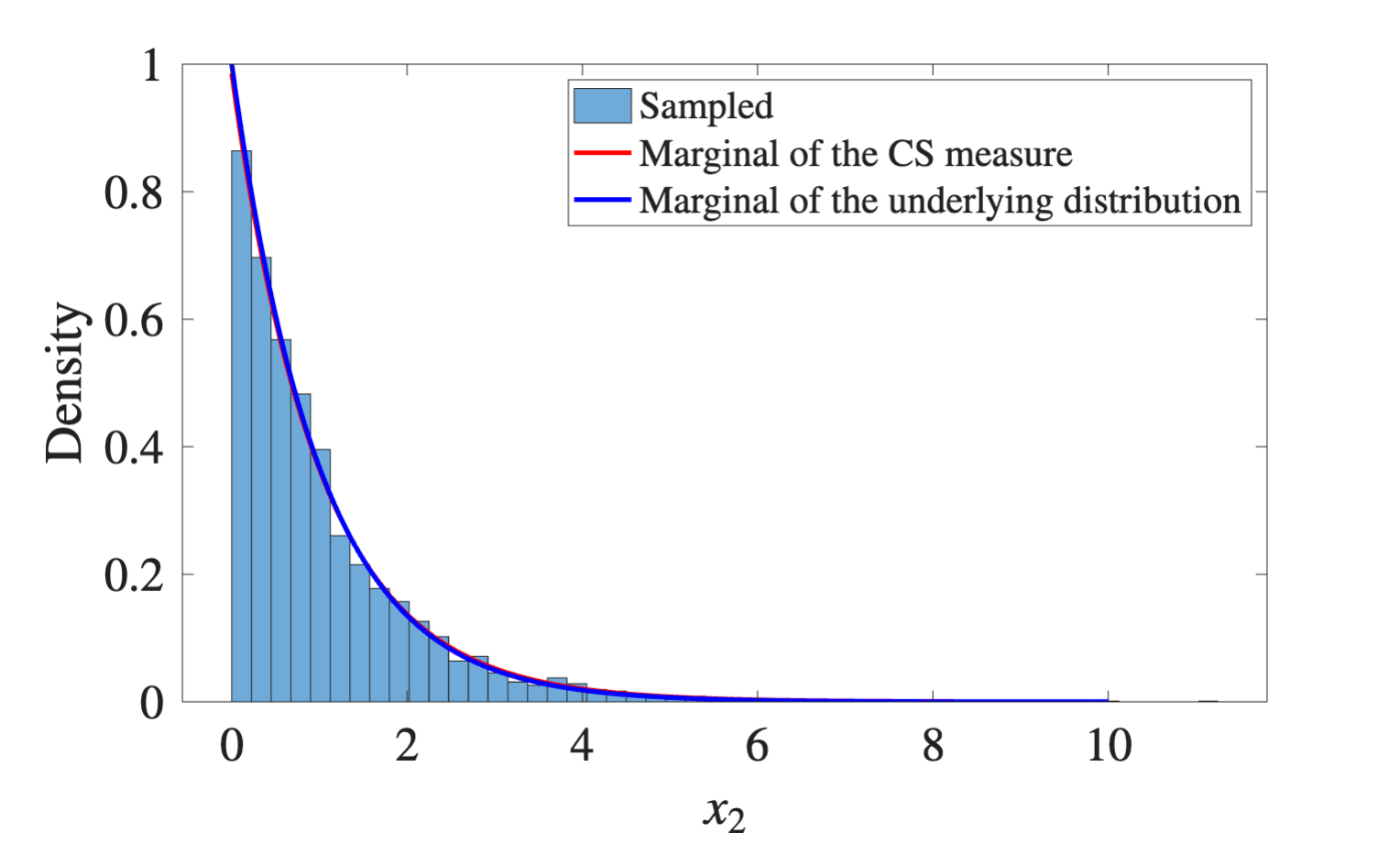}
        \caption{Marginal distribution of $x_2$.}
        \label{fig:mh2d_marginal_x2_exp}
    \end{subfigure}
    \caption{Diagnostics of the MH sampling algorithm with the CS measure \eqref{measurechris} with $\rho = \mathrm{Exp}(\lambda)$, $\lambda=1e-3$. The parameters are s = 100, $d$ = 2, $\sigma = 1$, $\sigma_1 = 1.5$, $\sigma_w=10^{-3}$, initial mean $x_0 = 1$, burn-in $B = 5000$, total samples $\texttt{total}= 5000$, thinning $T = 5$. }
\label{fig:2dmh_exp}
\end{figure}

Similar to the case of the normal distribution, we illustrate the MH sampler in the $d=1$ and $d=2$ cases, where where the underlying distribution is the multivariate exponential distribution $\rho = \mathrm{Exp}(\lambda)$ with the same parameter $\lambda > 0$ in each coordinate. As shown in Figures~\ref{fig:1dmh_exp} and \ref{fig:2dmh_exp}, the comparison between the empirical samples and the analytical density shows good agreement. Much as in the previous case, these figures also illustrate how the Christoffel function modifies the underlying exponential distribution.

\subsection{Synthetic examples}
\label{sec:numerics_syn}
\begin{figure}[htbp]
    \centering
    \begin{subfigure}[t]{0.45\textwidth}
        \centering
        \includegraphics[trim=0.8cm 0 0 0cm, clip, width=\linewidth]{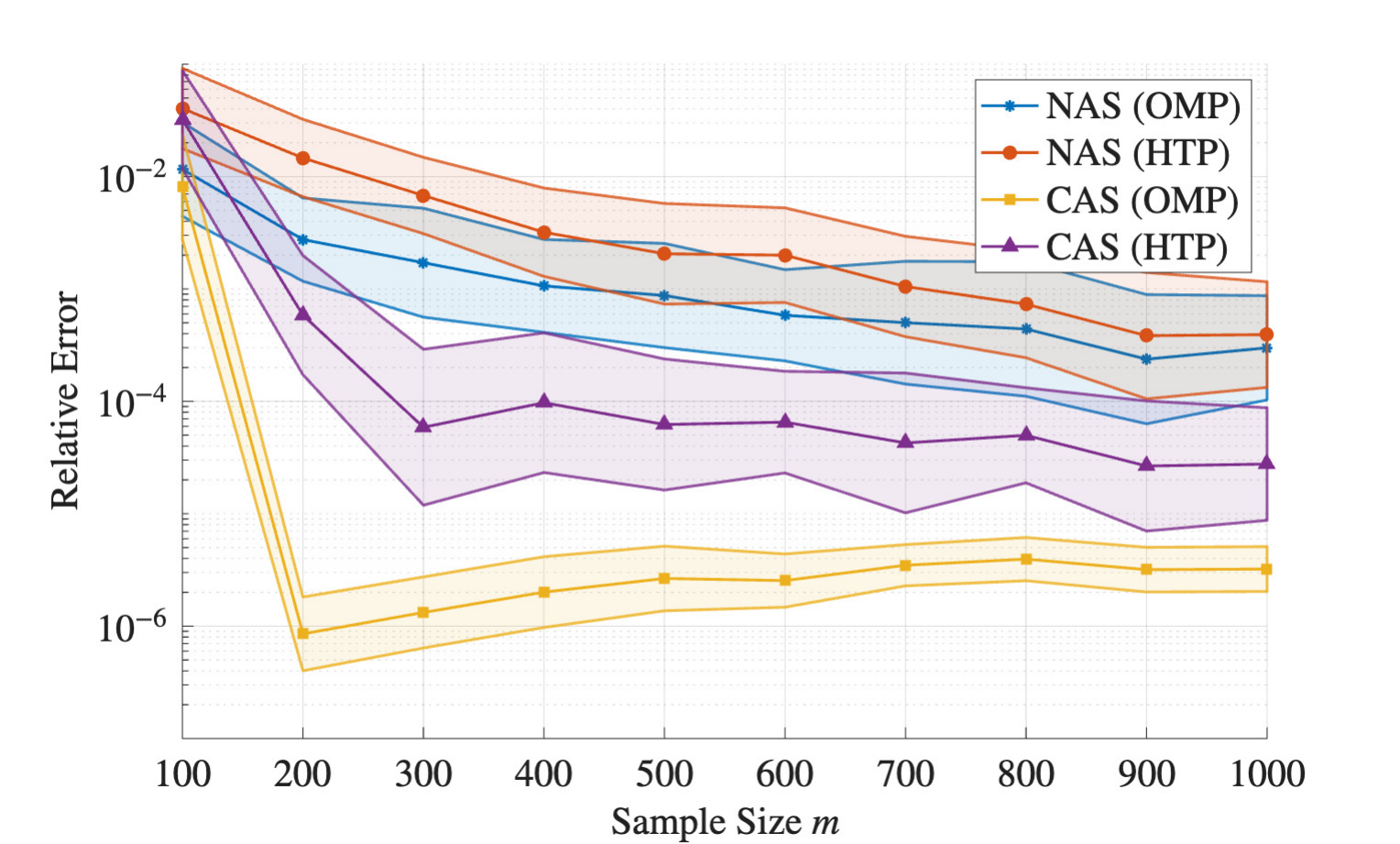}
        \caption{$f(x)=f_1(x)$, $d$ = 1; $N = 1000$}
        \label{fig:1d_exp}
    \end{subfigure}\hfill
    \begin{subfigure}[t]{0.45\textwidth}
        \centering
        \includegraphics[trim=0.8cm 0 0 0, clip, width=\linewidth]{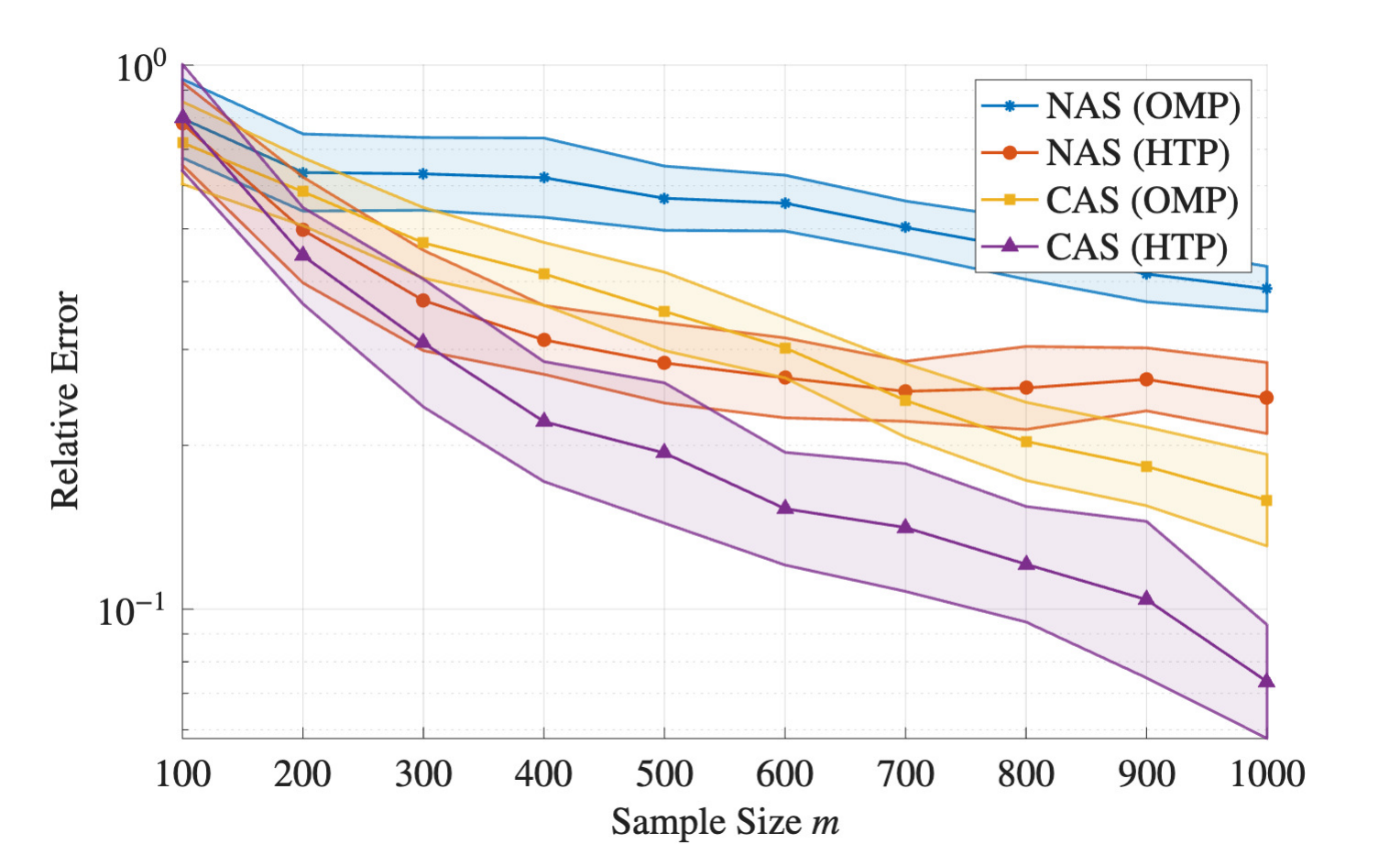}
        \caption{$f(x)=f_2(x)$, $d$ = 3; $N = 4000$}
        \label{fig:3dmh_sin}
    \end{subfigure}\hfill

     \begin{subfigure}[t]{0.45\textwidth}
        \centering
        \includegraphics[trim=0.8cm 0 0 0cm, clip, width=\linewidth]{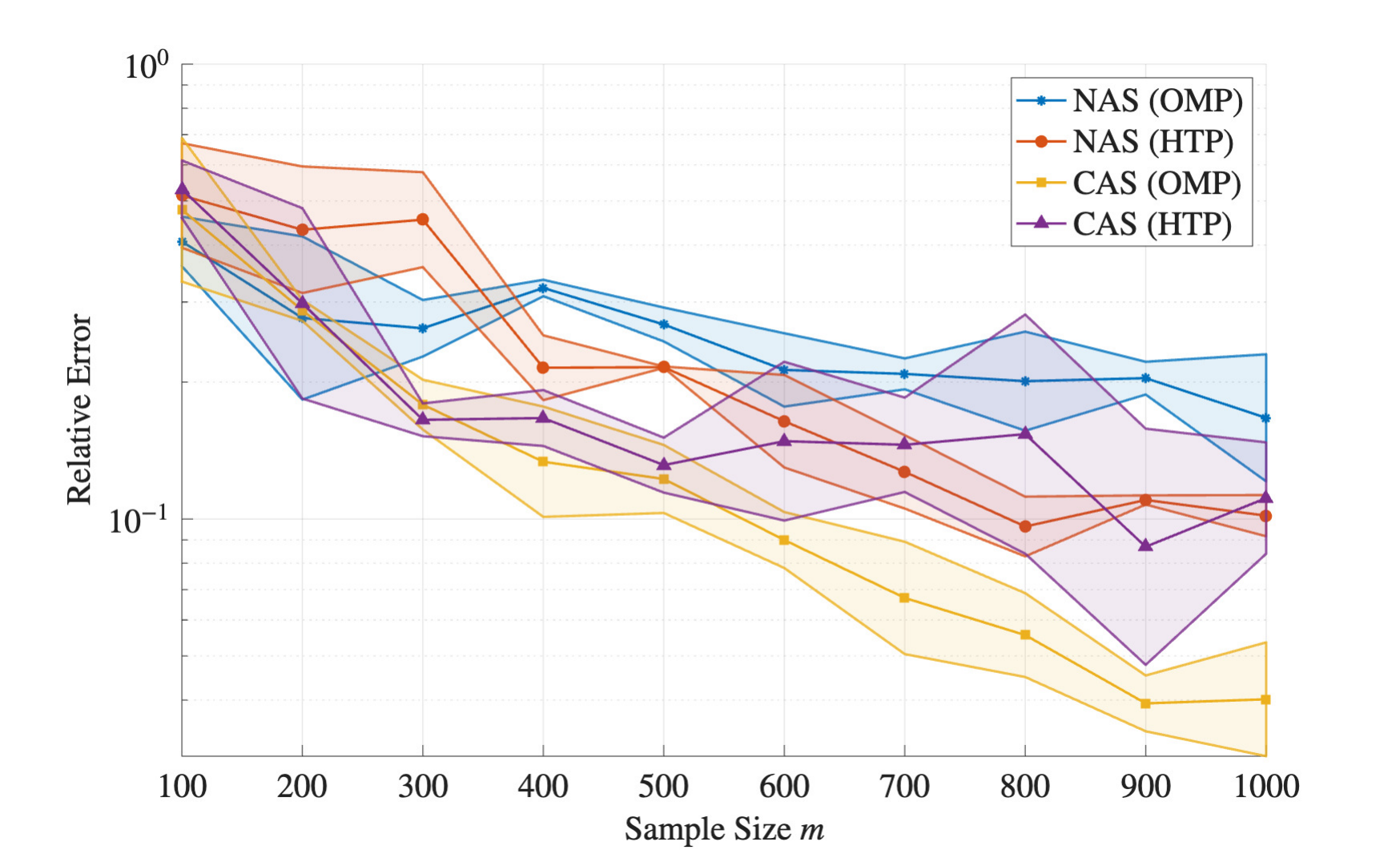}
        \caption{$f(x)=f_3(x)$, $d$ = 3; $N = 5000$}
        \label{fig:4dsumexp}
    \end{subfigure} \hfill
    \begin{subfigure}[t]{0.45\textwidth}
        \centering
        \includegraphics[trim=0.8cm 0 0 0cm, clip, width=\linewidth]{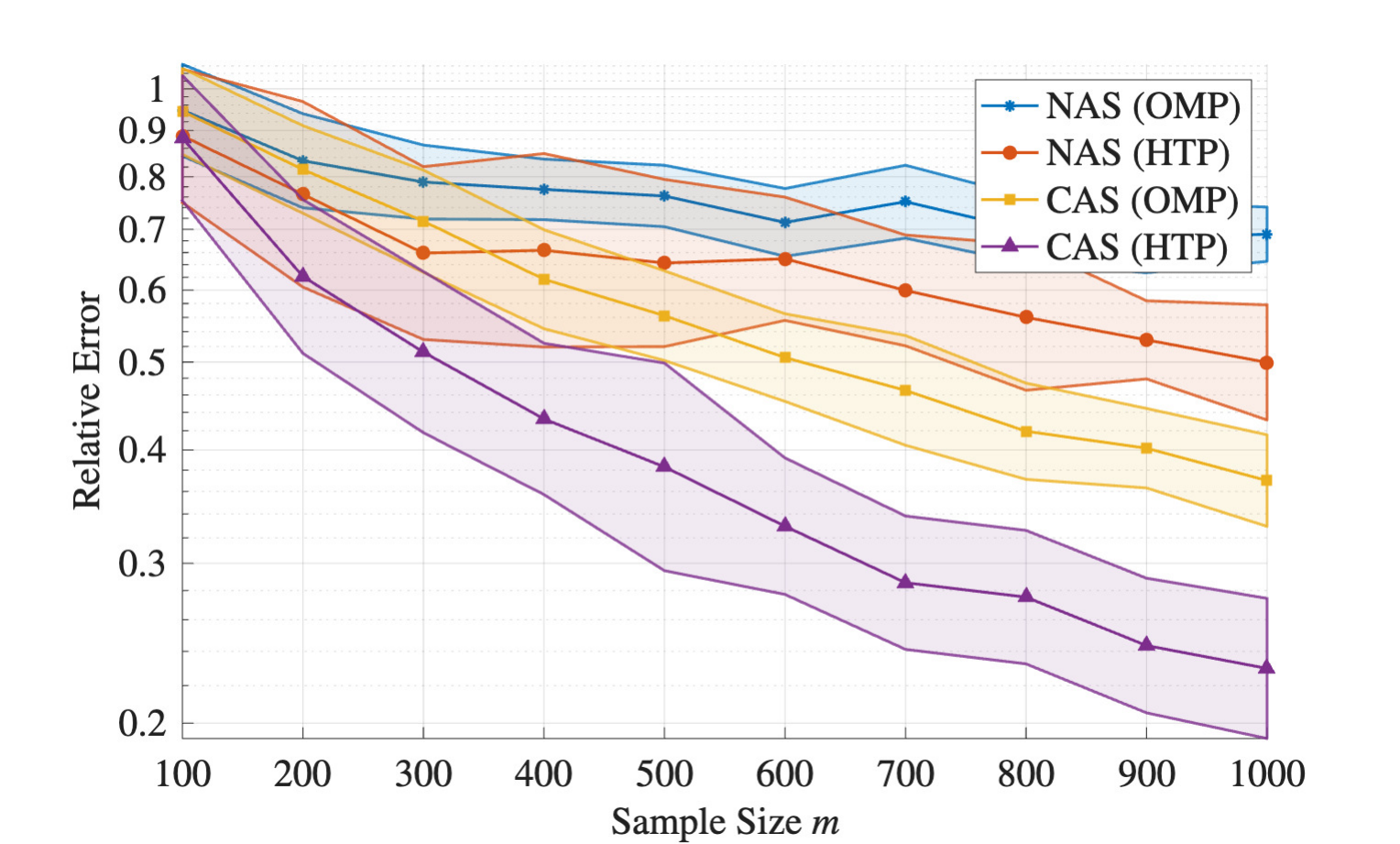}
        \caption{$f(x)=f_4(x)$, $d$ = 4; $N = 5000$}
        \label{fig:5dmh_sinsum}
    \end{subfigure}\hfill
    
    \begin{subfigure}[t]{0.45\textwidth}
        \centering
        \includegraphics[trim=0.8cm 0 0 0cm, clip, width=\linewidth]{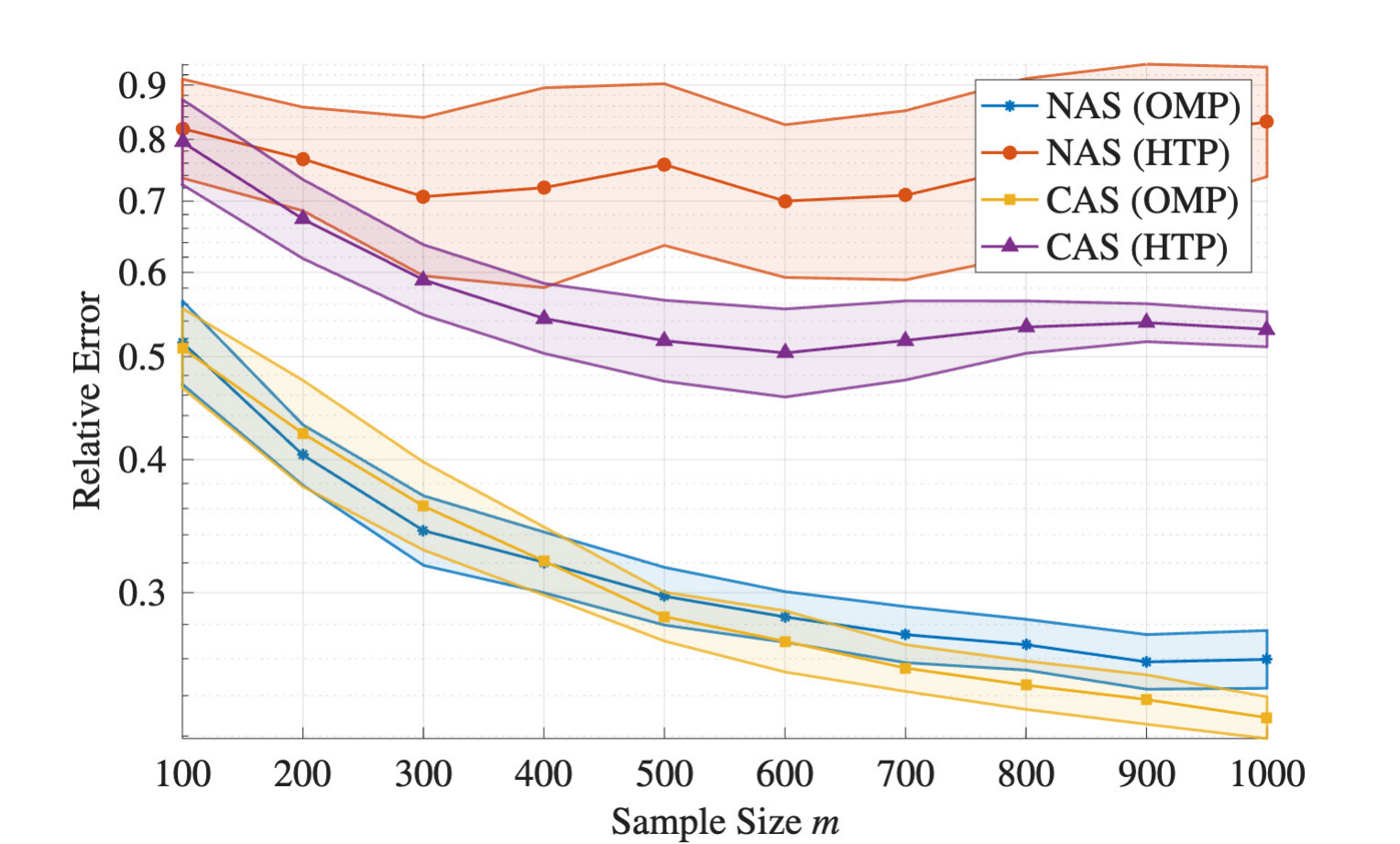}
        \caption{$f(x)=f_5(x)$, $d$ = 4, $N = 5000$}
        \label{fig:5dexpsum}
    \end{subfigure}\hfill
    \begin{subfigure}[t]{0.45\textwidth}
        \centering
        \includegraphics[width=\linewidth]{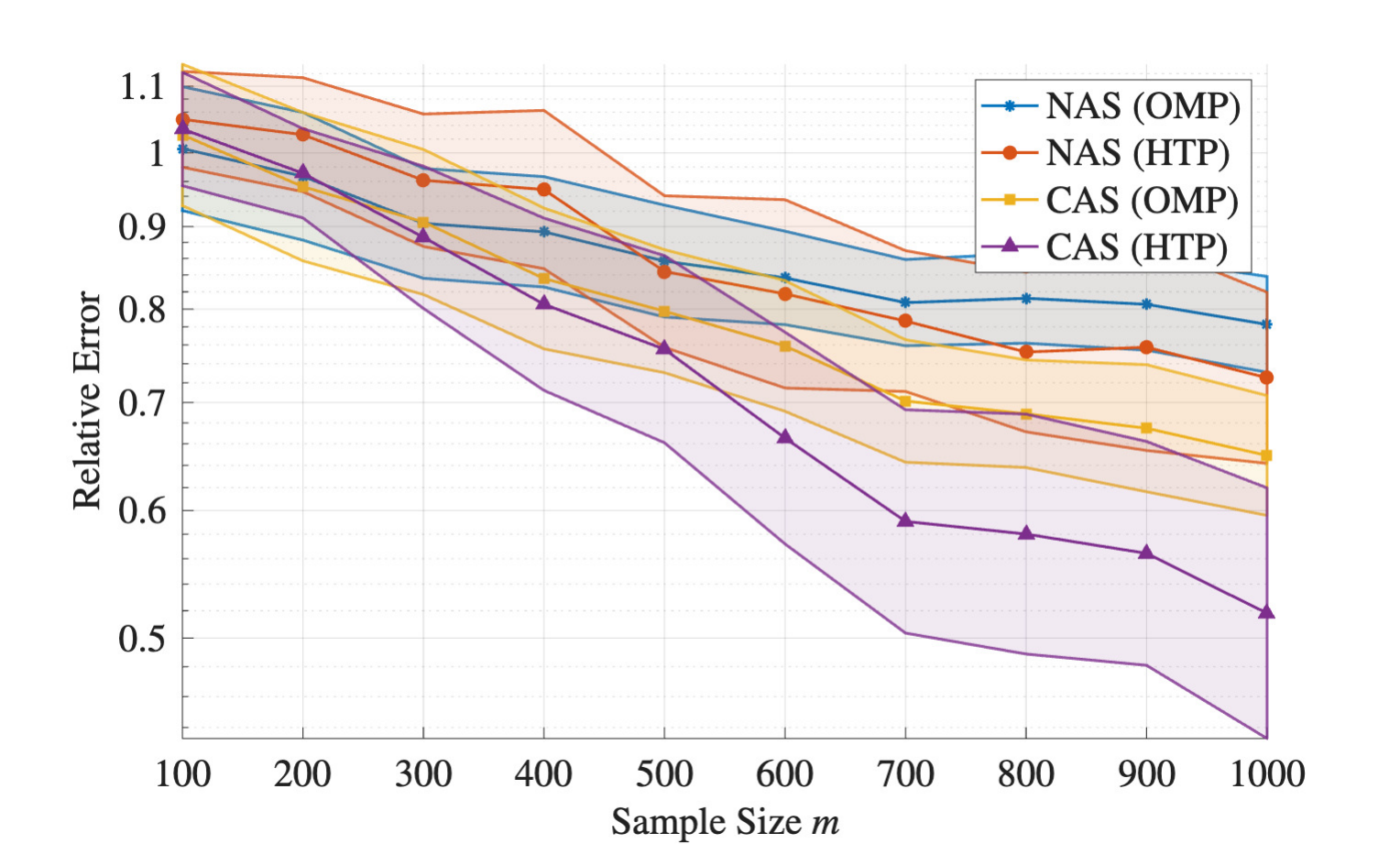}
        \caption{$f(x)=f_6(x)$, $d$ = 5; $N = 6000$}
        \label{fig:5dxsquaresum}
    \end{subfigure} \hfill
    \caption{CAS-SRFE and NAS-SRFE approximation errors versus $m$ for $\rho = \mathcal{N}(0,\sigma^2 I)$ with $\sigma = 1$. All experiments use sample sizes \texttt{m} = [100: 100: 1000]. Dimensions and number of random features are noted in each subfigure.}
    \label{fig:CASvsNAS}
\end{figure}

We now illustrate the performance of CAS-SRFE on a series of approximation tasks. In this section, we consider a suite of synthetic functions.
Let $\rho = \mathcal{N}(0,\sigma^2 I)$. In this case, we consider the following examples:
\begin{align*}
f_1(x) &= \exp(x), \quad d=1. \\
f_2(x) &= \sin(x_1) + 7 \sin^2(x_2) + 0.1\, x_3^4 x_1, \quad d=3. \\
f_3(x) &= \tfrac{1}{10} \sum_{j=1}^{d} \dfrac{\exp(-x_j^2)}{1+x_{j+1}^2}, \quad d=3.\\
f_4(x) &= \sum_{i=1}^d \Bigl(0.3 + \sin\!\left(\tfrac{16 x_i}{15} - 0.7\right) + \sin^2\!\left(\tfrac{16 x_i}{15} - 0.7\right)\Bigr),\quad d=4.\\
f_5(x) &= \sum_{i=1}^d \exp(-|x_i|),\quad d=5.\\
f_6(x) & = \left(\sum_{i=1}^d x_i\right)^2, \quad d=5.
\end{align*}
In Figure \ref{fig:CASvsNAS}, we report the error versus $m$ for both CAS and NAS, and for both the OMP and HTP sparse recovery algorithms. As shown, in the majority of cases, CAS leads to a consistent improvement over NAS. In many cases, the error is often over 10 times smaller than that of the corresponding NAS approximation. Note that the HTP variant often outperforms the OMP variant. However, this is not always the case, as can be seen in Figure~\ref{fig:CASvsNAS}(e). On the other hand, there are cases where CAS offers less benefit. This is especially prevalent as the dimension increases. However, this goes hand-in-hand with a worse approximation error. When $d = 1$, SRFE achieves higher accuracy, while when $d = 5$ it struggles to improve as the number of sample increases. Thus, the lack of improvement of adaptive sampling likely stems from the fact that RFMs do not provide sufficient expressibility for accurately approximating these higher-dimensional functions. In general, adaptive sampling will fail to provide any improvements if the underlying approximation scheme is not sufficiently accurate.

\begin{figure}[htbp]
    \begin{subfigure}[t]{0.5\textwidth}
        \centering
        \includegraphics[trim=0.8cm 0 0 0cm, clip, width=\linewidth]{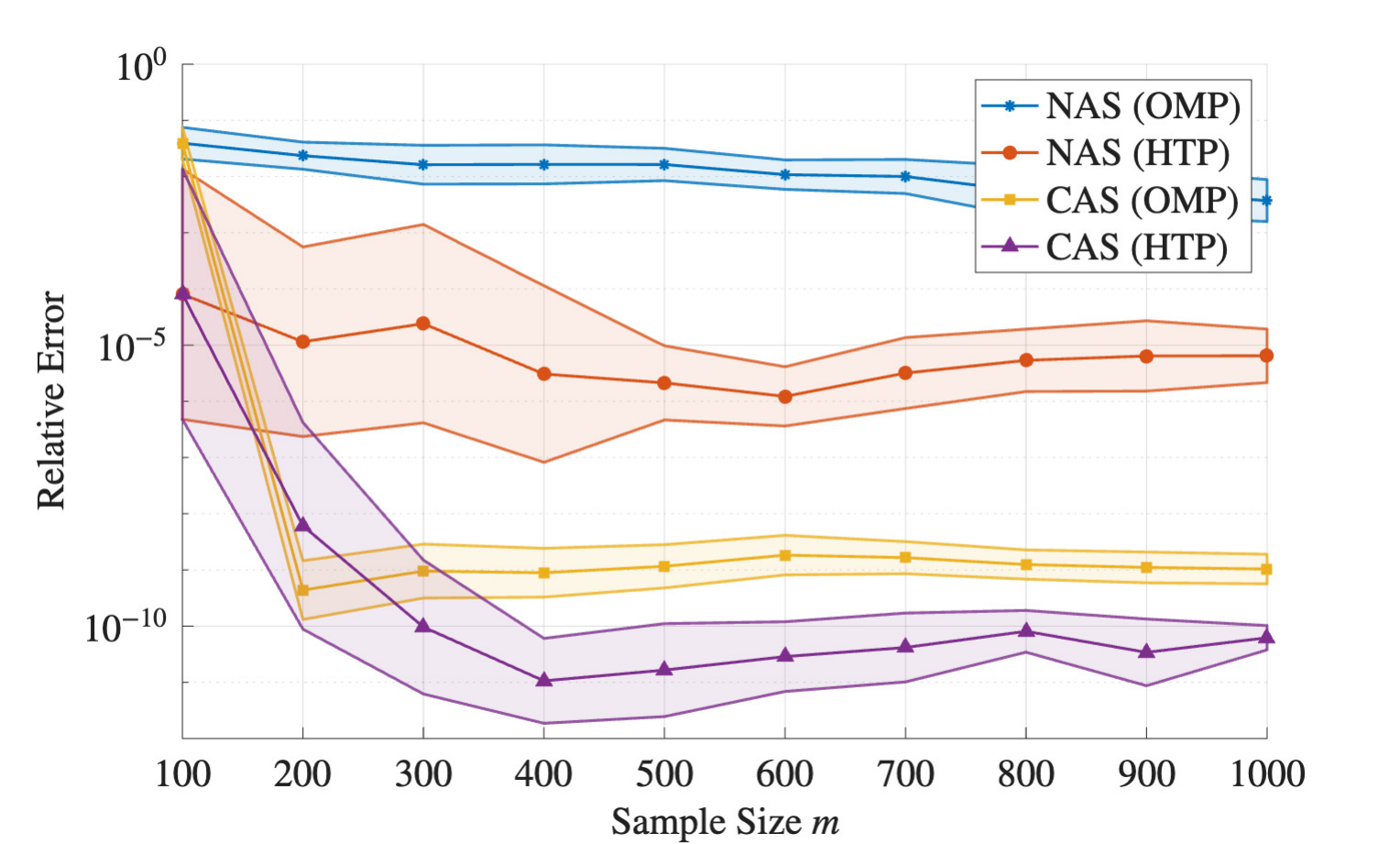}
        \caption{$f(x)=f_7(x)$, $d=1$, $N=1000$,  $\sigma_w=1$.}
        \label{fig:f7}
    \end{subfigure}\hfill
    \begin{subfigure}[t]{0.5\textwidth}
        \centering
        \includegraphics[width=\linewidth]{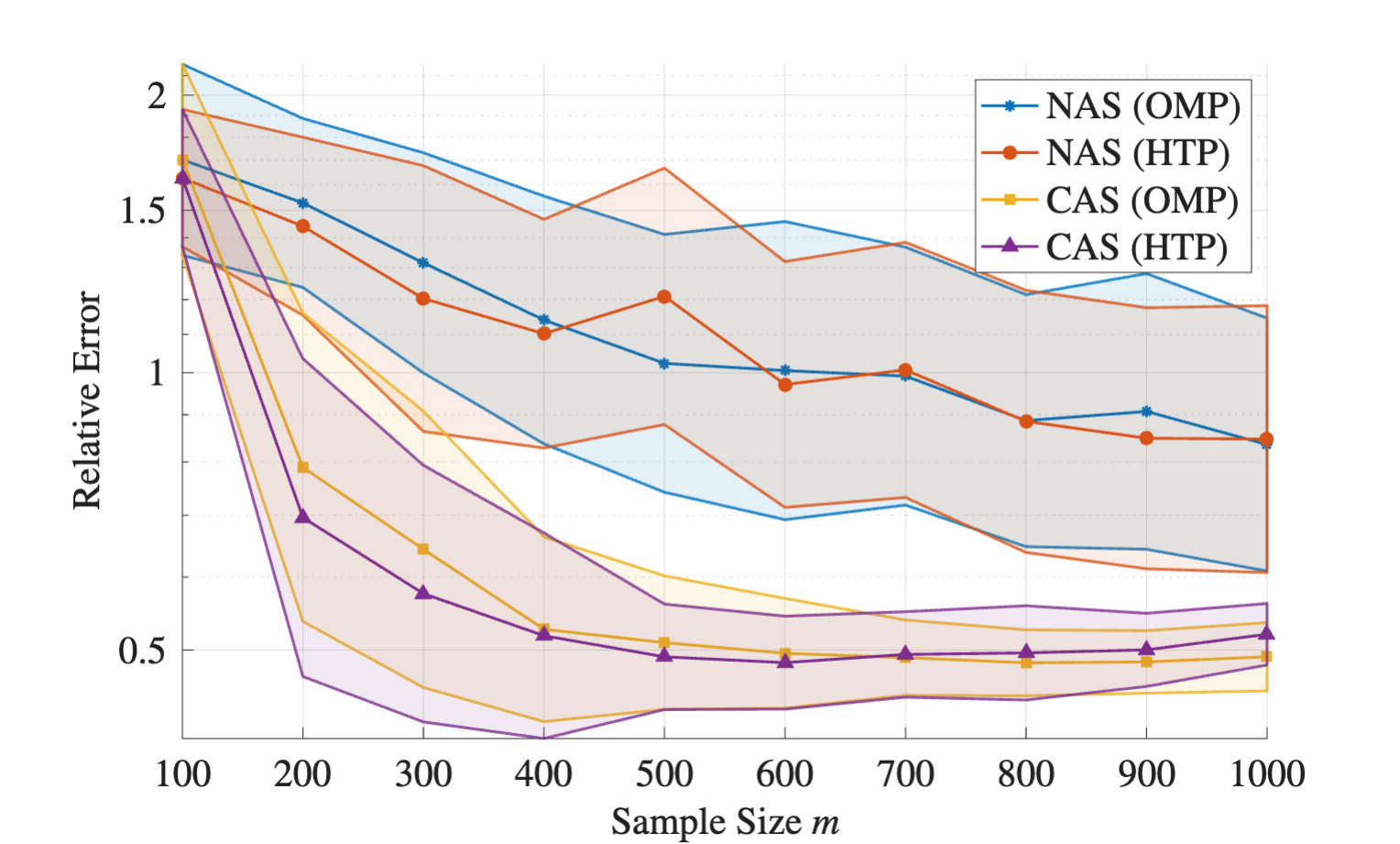}
        \caption{$f(x)=f_8(x)$, $d=2$, $N=3000$, $\sigma_w=1e-3.$}
        \label{fig:f8}
    \end{subfigure} \hfill
    \caption{CAS with MH sampling for samples with exponential underlying distribution. The solid line is the mean and the shaded region shows one standard deviation. All experiments use initial point $x_0$ = 1, sample sizes \texttt{m} = [100:100:1000]. Dimensions and number of random features are noted in each subfigure.}
    \label{fig:mh_comparison_exp}
\end{figure}

We also test for CAS-SRFE when the underlying distribution $\rho = \mathrm{Exp}(\lambda)$ is an exponential distribution. We consider the following examples:
\begin{align*}
f_7(x)=\sin(x), \quad d=1 \qquad
f_8(x)=\sin(x_2+x_2), \quad d=2. 
\end{align*}
As illustrated in Figure~\ref{fig:mh_comparison_exp}, CAS shows clear improvement over NAS. In both the one- and two-dimensional cases, CAS performs significantly better than NAS. Once more, we see a degradation in performance as the dimension increases.

\subsection{Parametric DEs}
\label{sec:parametric}
To investigate the effectiveness of CAS-SRFE in a more realistic setting, we apply it in the context in parametric modelling of physical systems, where the goal is to construct computationally efficient surrogate models. We consider several parametric Differential Equations (DEs), specifically, a surface reaction model, the Duffing oscillator and a damped harmonic oscillator. In all cases, the target function is a quantity-of-interest related to the parameters-to-solution map of a DE.

\subsubsection*{Surface Adsorption}
\label{sec:surface}
This example quantifies the uncertainty in the solution $\rho$ of the non-linear evolution equation. As discussed in \cite{hampton2018basis, hampton2015coherence}, this system models the surface coverage of certain chemical species. The governing equations are
\begin{equation}
    \begin{cases}
        \frac{d\rho}{dt}=\alpha(1-\rho)-\gamma \rho-k(1-\rho)^2\rho, \\
        \rho(t=0)=0.9,
    \end{cases}
\end{equation}
where $\kappa$ is the reaction constant, which is set to be $\kappa=10$. The adsorption coefficient $\alpha$ and desorption coefficient $\gamma$ are expressed as
\[\alpha=0.1+\mathrm{exp}(0.05 \Xi_1),\quad \gamma=0.001+0.01\mathrm{exp}(0.05 \Xi_2),\] where $\Xi_1, \Xi_2 \sim \mathcal{N}(0,1)$. The Quantity of Interest (QoI) is $\rho(t=4, \Xi_1, \Xi_2)$.
As illustrated in Figure~\ref{fig:surface}, CAS outperforms NAS. In this setting, only a relatively small number of random features is required, while the relative error achieved by CAS is approximately one order of magnitude smaller than that of NAS. 

\subsubsection*{Duffing Oscillator}
\label{sec:duffing}
Next we consider the Duffing oscillator, where $\rho(t)$ quantifies the displacement of a non-linear Duffing oscillator subject to an initial displacement under free vibration, as described in \cite{hampton2018basis}, governed by
\begin{equation}
\begin{cases}
\frac{d^2\rho}{dt^2}+2\chi \frac{d\rho}{dt}+\omega(\rho-\epsilon \rho^3)=0, \\
\rho(0)=1,\
\frac{d\rho}{dt}(0)=0.
\end{cases}
\label{eqn: duffingoscillator}
\end{equation}
Assuming the coefficients $\epsilon, \omega$ and $\chi$ are random and described by:
\begin{align*}
    \epsilon:=\mathrm{exp}(-0.1 \Xi_1^2), \quad \omega:=\mathrm{exp}(-0.1 \Xi_2^2), \quad \chi:=\mathrm{exp}(-0.1\Xi_3^2),
\end{align*}
where $\Xi_i \sim \mathcal{N}(0,1)$ for $i=1,2,3.$ As in \cite{hampton2018basis}, we consider the QoI $\rho(t=4 , \Xi_1,\Xi_2,\Xi_3)$. As shown in Figure~\ref{fig:duffing}, CAS consistently outperforms NAS with both OMP and HTP. This is expected as the estimation for solution of duffing model

\subsubsection*{Damped harmonic oscillator}
\label{sec:UHO}
\begin{figure}[htbp]
    \centering
    \begin{subfigure}[t]{0.48\textwidth}
        \centering
        \vspace{0pt} 
        \includegraphics[width=\linewidth]
        {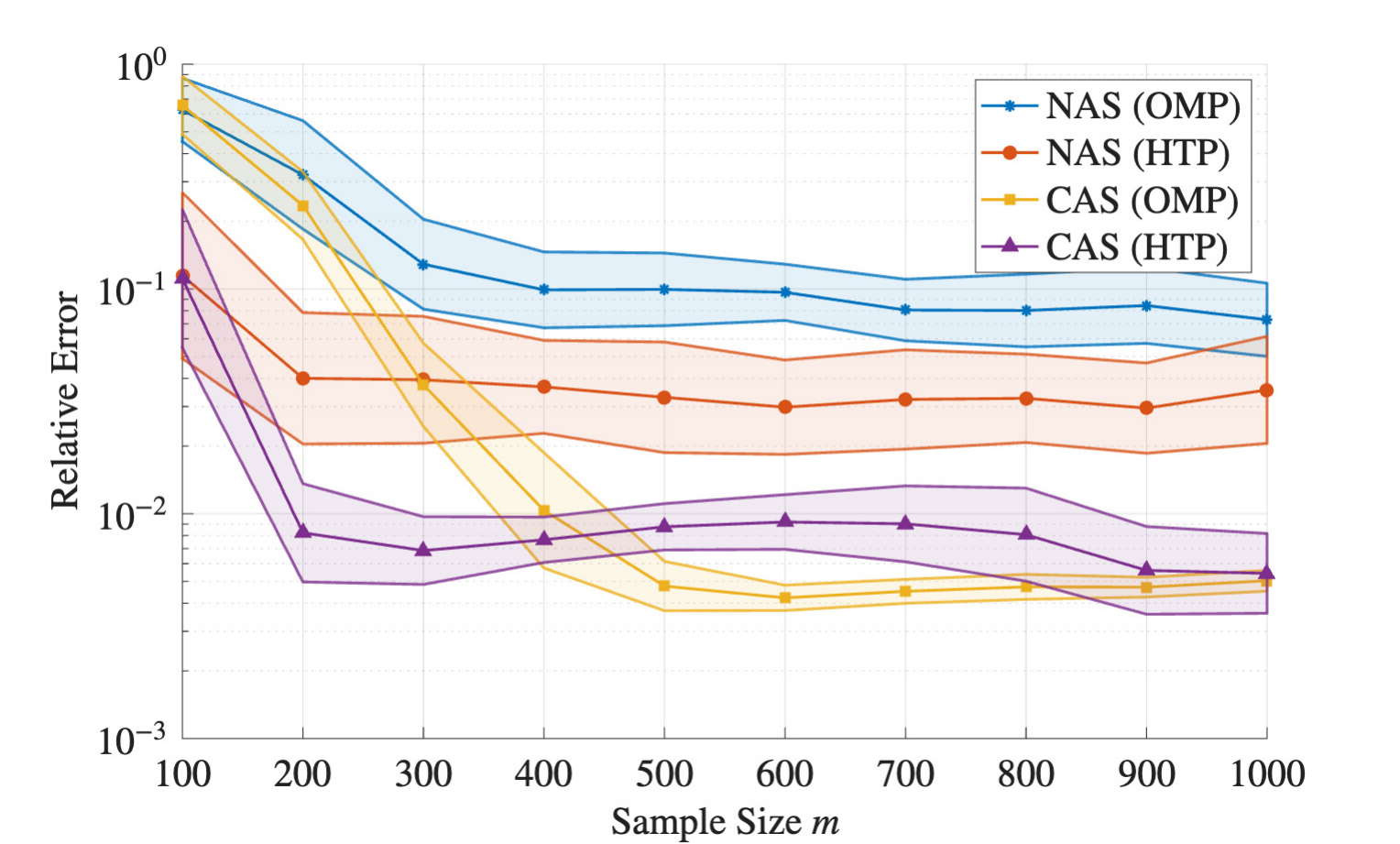}
        \caption{Surface adsorption model, $ N=2000. $}
        \label{fig:surface}
    \end{subfigure}
    \hfill 
    \begin{subfigure}[t]{0.48\textwidth}
        \centering
        \vspace{0pt}
        \includegraphics[width=\linewidth]{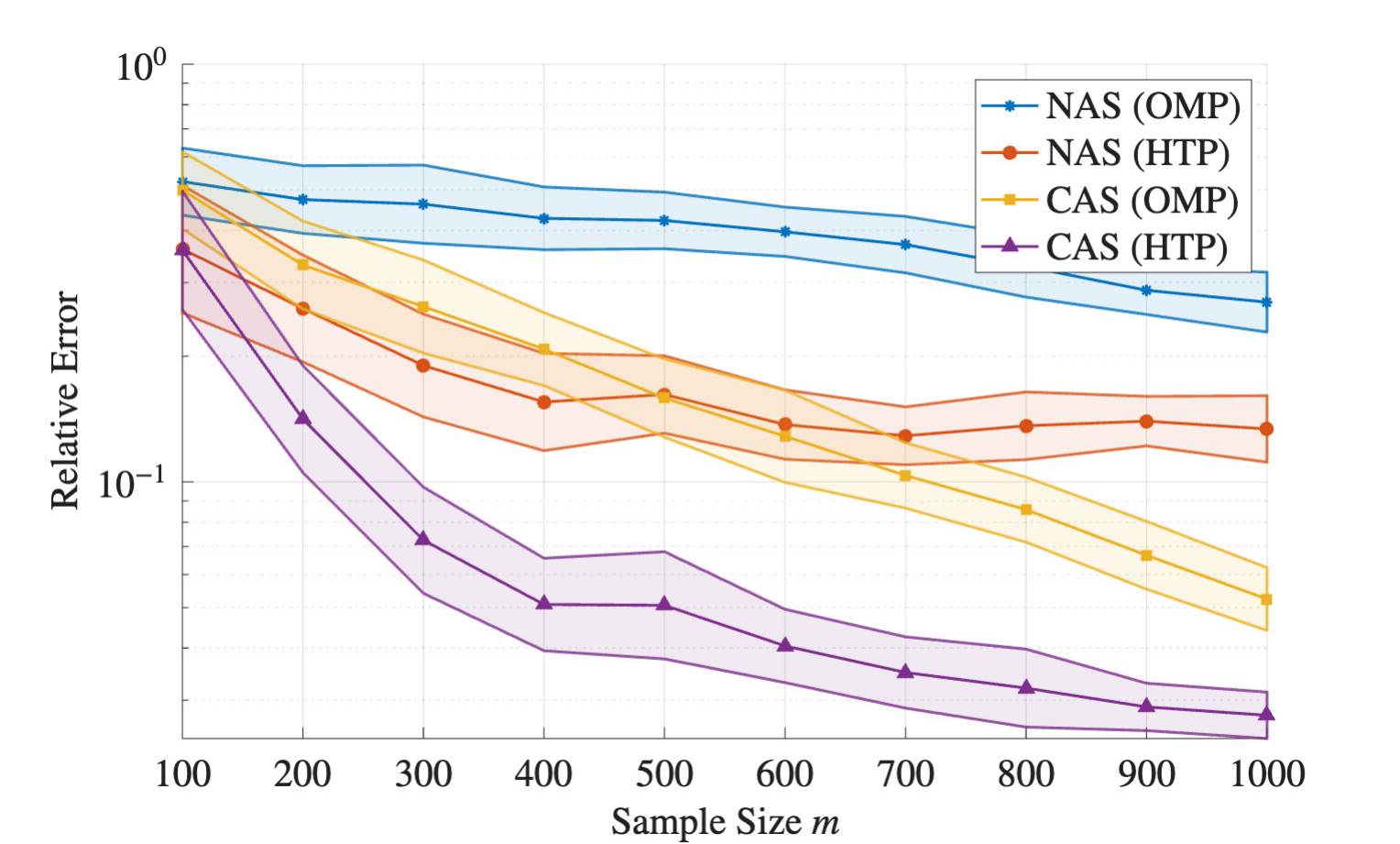}
        \caption{Duffing oscillator model, $N=4000$.}
        \label{fig:duffing}
    \end{subfigure}
    \hfill 
    \begin{subfigure}[t]{0.48\textwidth}
        \centering
        \vspace{0pt}
        \includegraphics[width=\linewidth]{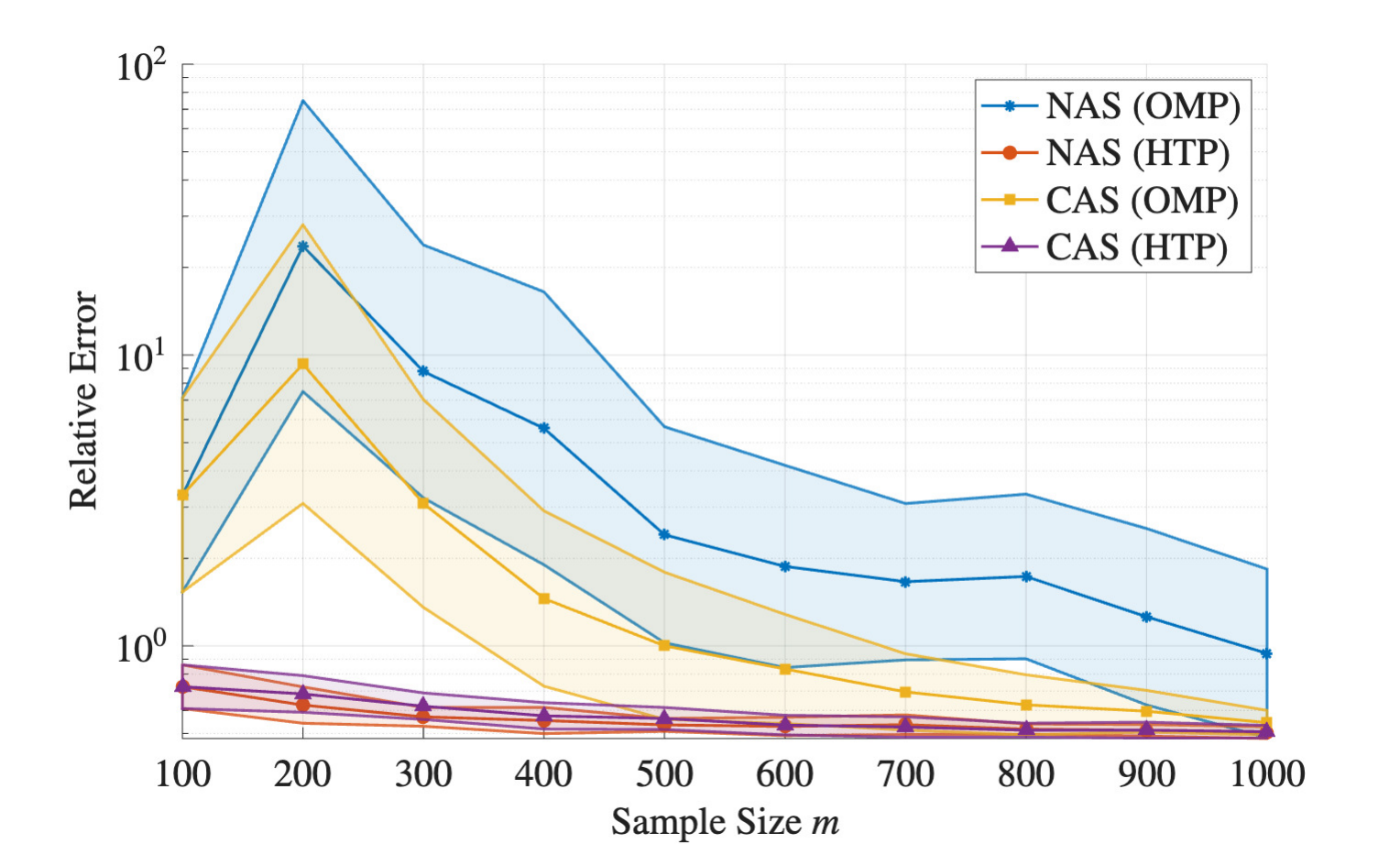}

        \caption{Underdamped harmonic oscillator, $x_0=1$, $\sigma_w = [1,1,10^{-3},10^{-3},10^{-3},1]$, $N=8000$. }
        \label{fig:harmo_fig}
    \end{subfigure}

    \caption{Comparison of CAS-SRFE and NAS across three different systems: (a) underdamped harmonic oscillator, (b) Duffing oscillator, and (c) surface adsorption model. All experiments use sample sizes \texttt{m} = [100:100:1000].}
    \label{fig:three_systems_comparison}
    
\end{figure}

We investigate a system of an underdamped harmonic oscillator. The system is subject to external forcing with six unknown parameters: damping coefficient $\gamma$, spring constant $k$, forcing amplitude $g$, frequency $\omega$, and the initial conditions $u_0$ and $u_1$. 
Defining \[
\Xi = (\gamma, k, g, \omega, u_0, u_1)
   = (\Xi_1, \Xi_2, \Xi_3, \Xi_4, \Xi_5, \Xi_6).
\] 
As such, we express the governing equations as following:
\begin{equation}
\begin{cases}
\frac{d^2 u}{dt^2}(t,\Xi) + \gamma \frac{d u}{dt} + k u
    = g \cos(\omega t), \\
u(0,\Xi) = u_0,\ 
\frac{du}{dt}(0,\Xi) = u_1.
\end{cases}
\label{eqn: harmo}
\end{equation}
To guarantee that the system is underdamped, one needs $\gamma^2-4k<0$. For implementation, modifying \cite{adcock2018compressed}, we set $\Xi^{(1)}=0.1, \quad \Xi^{(2)}=0.04$. Let $\Xi^{(i)}\sim k_i+\operatorname{Exp}(\lambda_i)$ for $i=3,4,5$, $\Xi^{(6)}\sim \mathcal{N}(0,0.5)$. Fix the parameters $(k_3, \lambda_3)=(0.08, 1), (k_4, \lambda_3)=(0.8. 1), (k_5, \lambda_5)=(0.45, 1)$. As in \cite{adcock2018compressed}, we consider the QoI $u(t=20, \Xi)$. In Figure~\ref{fig:harmo_fig}, we present the results for estimation. 

Since this example is rather complicated in terms of distribution and has higher dimensions, the number of random features required is much higher. Both NAS and CAS fail to  achieve high accuracy, indicating that SRFEs are unable to approximate the underlying function well. However, CAS still achieves a modest improvement.

\section{Conclusions and future work}
\label{sec:conclusions}

By combining an adaptive sampling process guided by CAS with the sparsity-based representation model SRFE, we introduced SRFE-CAS for efficient function approximation in data-scarce settings. We introduced a novel method to draw samples from the relevant CS measures, and then conducted a series of experiments to demonstrate the performance of SRFE-CAS. Our experiments show a consistent improvement of NAS, However, in higher-dimensional experiments, both CAS and NAS fail to provide accurate approximations. This is likely due to the increase in dimensionality, which makes it more difficult to tune the underlying SRFE procedure to achieve high accuracy. 

There are several avenues for future work. First, our framework requires explicit computation of a Gram matrix, which limits its applicability to simple domains and RFMs. However, recent work \cite{herremans2025refinement} proposes an approach that enables computation of the CS measure without explicit orthogonalization, starting from a non-orthogonal basis and using an iterative refinement strategy. Combining this with CAS-SRFE is an interesting option for future work. Second, CAS–SRFE provides a general framework that is not restricted to the representation model introduced in this paper. Alternative sparse recovery algorithms may be used in place of HTP or OMP, and other nonlinear approximation schemes beyond RFM-based models may also be considered. For instance, the algorithm could be further enhanced by exploring alternative representation frameworks. For example, \cite{hashemi2023generalization} introduces Sparse Random Feature Expansion with Sparse Features (SRFE-S), which enforces sparsity in the feature weights by sampling $q$-sparse random features—an approach not explored here. In addition, RFMs have recently been applied in various scientific computing contexts, such as numerical PDE solvers \cite{chen2022bridging, liao2026solving} and operator learning \cite{nelsen2024operator}. It would be interesting to explore CAS-SRFE for these problems. More broadly, CAS-SRFE highlights the importance of well-designed sampling strategies and representation models for approximation tasks with data scarcity, which is a common issue for many scientific problems. This work provides a promising direction for building tractable, scalable, and efficient approximation methods in scientific computing.

\bibliographystyle{siamplain}
\bibliography{references.bib}
\end{document}

%% file: ex_shared.tex
\usepackage{lipsum}
\usepackage{amsfonts}
\usepackage{graphicx}
\usepackage{epstopdf}
\usepackage{algorithmic}
\ifpdf
  \DeclareGraphicsExtensions{.eps,.pdf,.png,.jpg}
\else
  \DeclareGraphicsExtensions{.eps}
\fi

\newsiamremark{remark}{Remark}

\newsiamremark{hypothesis}{Hypothesis}
\crefname{hypothesis}{Hypothesis}{Hypotheses}
\newsiamthm{claim}{Claim}
\newsiamremark{fact}{Fact}
\crefname{fact}{Fact}{Facts}

\headers{Christoffel Adaptive Sampling for Sparse Random Feature Model}{}

\title{Christoffel Adaptive Sampling for Sparse Random Feature Expansions
\thanks{Submitted to the editors \today.
\funding{
 BA acknowledges support from the Natural Sciences and Engineering Research Council of Canada (NSERC) through grant RGPIN/2470-2021. BA \& XW acknowledge the support of FRQ (Fonds de recherche du Qu\'ebec) -- Nature et Technologies through grant 359708.}
}}

\author{
  Ben Adcock\thanks{Department of Mathematics, Simon Fraser University 
    (\email{ben\_adcock@sfu.ca}).}
  \and
  Khiem Can\thanks{Department of Mathematics, Simon Fraser University 
    (\email{dinh\_khiem\_can@sfu.ca}).}
    \and
      Xuemeng Wang\thanks{Department of Mathematics, Simon Fraser University 
    (\email{xwa265@sfu.ca}).}
}

